\documentclass[11pt]{article}
\usepackage{amsthm,amsmath,amsfonts}
\usepackage{latexsym}
\usepackage[T1]{fontenc}
\usepackage{a4wide}
\usepackage{times}
\usepackage{array}
\usepackage{hhline}
\numberwithin{equation}{section}

\input epsf

\theoremstyle{plain}
\newtheorem{theorem}{Theorem}[section]
\newtheorem{proposition}{Proposition}[section]
\newtheorem{lemma}{Lemma}[section]
\newtheorem{corollary}{Corollary}[section]
\theoremstyle{definition}
\newtheorem{definition}{Definition}[section]

\theoremstyle{remark}

\def\al{\alpha}
\def\ga{\gamma}
\def\Ga{\Gamma}
\def\ka{\kappa}
\def\Om{\Omega} \def\om{\omega}
\def\be{\beta}
\def\de{\delta}
\def\ep{\epsilon}
\def\vfi{\varphi}
\def\bfi{\bar\varphi}
\def\La{\Lambda} \def\bLa{\bar\Lambda}
\def\la{\lambda}
\def\De{\Delta}
\def\si{\sigma}

\newcommand{\bR}{\mathbb R}
\newcommand{\JR}{{\rm J}\mbox{-}{\rm R}}
\def\bn{{\bf n}}

\def\te{{\bf e}}
\def\TV{\hbox{TV}\,}
\def\cA{\mathcal A}
\def\cD{\mathcal D}
\def\cJ{\mathcal J}

\def\cM{\mathcal M}

\begin{document}
\title{A caricature of a singular curvature flow in the plane}
                                                                           
%%%%%%%%%%%
% Authors %
%%%%%%%%%%%
\author{Piotr B. Mucha and  Piotr Rybka\\
Institute of Applied Mathematics and Mechanics,
Warsaw University\\
ul. Banacha 2, 02-097 Warszawa, Poland\\
E-mail: p.mucha@mimuw.edu.pl, p.rybka@mimuw.edu.pl\\
Corresponding author: Piotr Rybka, p.rybka@mimuw.edu.pl, fax: +48 22 55 44 300
}
\maketitle

\date{}

\bigskip\noindent{\bf Abstract.}  We study a singular parabolic
equation of the total variation type in one dimension. The problem is
a simplification of the singular curvature flow. We
show existence and uniqueness of weak solutions. We also prove
existence of weak solutions to the semi-discretization of the problem
as well as convergence of the approximating sequences. The
semi-discretization shows that facets must form. For a class of initial
data we are able to study in details the facet formation and
interactions and their asymptotic behavior. We notice that our
qualitative results may be interpreted with the help of a special composition of
multivalued operators. 


\bigskip\noindent{\bf keywords:} singular parabolic equations, 
singular curvature flow, monotone operators, facet formation, facet interaction
%\bigskip

%%% --- s 1 ---

\section{Introduction}

Many free boundary problems involving the Gibbs-Thomson relation may be
considered as a driven weighted mean curvature flow  coupled
through the forcing term to a diffusion equation (see \cite{Chen-Re},
\cite{Radkevich}, \cite{luckhaus} \cite{al-wang}). We have a
considerable body of literature concerning this problem
for the Euclidean curvature of the interface, including the question of
precise regularity of solutions treated by Escher, Pr\"uss, Simonett and Mucha,
see \cite{eps}, \cite{epi}, \cite{mu1}. On the other hand,
less is known if the
curvature appearing in the Gibbs-Thomson relation is singular, see
e.g. \cite{rybka}.  This line of
research has been initiated by Taylor, \cite{taylor}, and independently by
Gurtin, \cite{gurtin}. However, just solvability of equations  of  the
singular curvature flow is interesting. Existence of the flow was
obtained by Bellettini, 
Novaga, Paolini \cite{bnp1}, \cite{bnp2}  and by Chambolle
\cite{chambole}. Driven singular curvature flow was studied by M.-H.Giga,
Y.Giga and  Rybka, see \cite{gigi}, \cite{gr1}, \cite{gr2}.

In fact, the existence and properties of solutions to
the  singular weighted mean  curvature flow
\begin{equation}\label{rn0}
 V = \ka \qquad \hbox{on \ \ }\Ga(t),
\end{equation}
are interesting in itself even in the plane and without forcing,
especially
when the anisotropy function (also called `energy density function')
is singular, i.e. just convex. Here,  $\Gamma(t)$ is the unknown curve and
$\kappa$ denotes the weighted mean  curvature related to the underlying
anisotropy function and $V$
is the velocity of surface $\Gamma(t)$.
Our ultimate  goal would be to study existence and behavior of solutions to
(\ref{rn0}).

In its full generality problem (\ref{rn0}) for an arbitrary initial curve is
rather difficult. One source of difficulties is
the geometry of the system, it is already present in the two-dimensional
setting. Here, we want to concentrate only on the purely
analytical difficulties appearing in (\ref{rn0}). This is why we  will
restrict our attention to a simplified equation, which retains the singular
character of the original problem.
%One simplification which can be found in
%the literature was to consider curves as graphs over intervals.

Here is our
postulated equation
\begin{eqnarray}\label{rn2}
%\begin{equation}
&&\La_t = %\frac d{ds}
\frac{\partial}{\partial s}
\frac d{d\phi}J%I_\theta
(s+\La_s)\quad \hbox{in }S\times(0,T),\nonumber\\%quad
&&\La(s,0) =
\La_0(s)\quad \mbox{on } S,\\%qquad
&&\La(2\pi,t) = \La(0,t), \quad t\ge 0,\nonumber
\end{eqnarray}%equation}
here $S$ is the unit circle parameterized by interval $[0,2\pi)$ and $\Lambda$
is the sought function. Compared with (\ref{rn0}) our new system has one
analytical advantage. Namely, the domain of definition of $\La(\cdot,t)$ is
independent of time.

We present a justification of this equation in the Appendix. Here, we explain
our notation. The variable $s$ plays the role of the arclength parameter, 
the subscript $s$ denotes the differentiation with respect to $s$. We
frequently refer to $\vfi=\La_s +s$ as the angle between the $x_1$ axis and
the outer normal to the curve. Such an interpretation helps drawing pictures,
but the relation to the actual angle is rather loose.

We  make a specific choice of $J$ corresponding to the surface energy
 density functions. We want to study a situation which is already very
 singular yet tractable. In many instances of a great physical interest an
 anisotropy appears, which is merely convex, not even strictly convex
 (understood in a proper sense). As a result, we
choose $J$, which is convex and piecewise linear. This is an independent
 source of difficulties.  In order to avoid further technical troubles we
 will choose $J$ corresponding to the situation where that curve minimizing
 the surface energy (which is the Wulff shape of the anisotropy function) is a
 square. We must stress again that the correspondence is at the level of
 ideas, because  (\ref{rn2}) is {\it not} a curvature flow, but its
 caricature. However, the obtained behavior of solutions to (\ref{rn2})
 is almost the same
 as for the equation (\ref{rn0}) with the anisotropy function
 corresponding to a square, \cite{chambole}.

Thus, we pick $J$ which suffers jumps of equal height $\frac\pi2$ at the
equi-spaced angles
\begin{equation}\label{angle}
{\cal A}=\left\{\al_k= -\frac{3\pi}4+k\De\al :\quad k=0,1,2,3,\qquad
\mbox{with } \De\al=\frac\pi2 
\right\}.
\end{equation}
Specifically, we put
\begin{equation}\label{jot}
J(\vfi) =
\frac\pi4\left(
|\vfi -\frac{3\pi}4 |+ |\vfi-\frac\pi4 |+
|\vfi +\frac\pi4 |+ |\vfi +\frac{3\pi}4 |
\right).
\end{equation}
Since $\La$ is defined over the unit circle its graph over $S$ is a closed
curve. The meaning of the spacing between $\al_k$'s can be explained by
looking at the equation
$$
\frac \partial{ \partial s}\frac d{d\vfi} J(s+\La_s) =1,
$$
considered in \cite{MuchaRybka} -- see subsection 3.2, too. Roughly speaking,
the spacing between $\al_k$ and $\al_{k+1}$ corresponds to the length of
facets having the normal vector $\bn$ with the normal angle $\al_k$. The size
of the jump of 
$\frac d{d\vfi} J(s+\La_s)$ corresponds to the angle between the normals to the
curve, which is a solution to the above equation, at a corner.

The chosen  anisotropy function (\ref{jot})  is nowhere regular, hence we can expect
nonstandard effects  requiring  new analytical tools. This has been observed
by researchers working on the total variation flow, whose simplification is
\begin{equation}\label{tvf}
u_t-\delta_0(u_x)u_{xx}=0
\end{equation}
augmented with initial and boundary data. Here, $\delta_a$ is the Dirac measure
concentrated at $a$.

We noticed so far two main types of motivation to study the total variation
flow,
\begin{equation}\label{tvfn}
u_t-\hbox{div}\,\left(\frac{\nabla u}{|\nabla u|}\right)=0.
\end{equation}
The first one is the image denoising and reconstruction introduced by
Rudin and Osher, \cite{ruso1}, \cite{ruso2}. The second one is evolution of
the facets of crystals. The bulk of the papers (see \cite{andreu-bcm1},
\cite{andreu-bcm2}, \cite{bcn}, \cite{gigakoba}, \cite{gigikoba},
\cite{acm-book}, \cite{moll}) uses the theory of nonlinear semigroups
to establish existence. The last paper is particularly interesting,
because it deals with the anisotropic  total variation
flow. Moreover, the notion of entropy solutions was
introduced to deal with uniqueness of the total variation flow (see
\cite{andreu-bcm1}, \cite{bcn}).
The tools of convex analysis were useful to make sense out of (\ref{tvf}).
The authors, mentioned above, paid special attention to piecewise
constant initial data and they were interested in the asymptotic
behavior, in particular the asymptotic shape was identified.
M.-H.Giga, Y.Giga and R. Kobayashi, \cite{gigakoba},
\cite{gigikoba}, also calculated the speed of flat facets.
No matter what is the approach, it is apparent that the most important
information  is located in sets $\{u_x=0\}$, where the singular dissipation
starts to play a role and where the classical multivalued theory of
function loses the meaning.

%\bigskip
Our approach differs in many aspects. We prove existence by a regularizing
procedure and passing to the limit with the regularizing parameter,
this approach was used, e.g. by Feng and Prohl, see \cite{fengprohl}.
The main difficulty is associated with studying the limit of the non-linear
terms. We present a more detailed analysis of regularity of solutions
permitting us to call them `almost classical'. For generic data, our solutions
are twice differentiable with respect to $s$, except a finite number of points
(for fixed time).
This will be explained in detail below. We mention here that we use the tools
of the convex analysis, in particular we rely on the fact that for a
convex function the subdifferential is well-defined everywhere.
However, the classical theory of multivalued functions is not sufficient. We
have to introduce a
new definition of the composition of two multivalued functions to describe  the
meaning and qualitative properties of solutions to system (\ref{rn2}) as well
a class of the $\JR$ functions, where regularity is described from the point
of view of the properties of the function $J$. In our opinion the
results we prove contribute  to better understanding parabolic systems with
measure coefficients.

Our technique requires a new look at the regularity of functions. We
will generalize the meaning of the convexity defining
a class of $J$-regular functions  preserving some important properties of the
convexity. Our main qualitative result says that any sufficiently regular
initial curve evolving according to system (\ref{rn2}), will eventually reach
a minimal solution, which is called  the asymptotic profile in the area of the
total variation flow. The geometric interpretation
is that the solution reaches its asymptotic shape, i.e. the square in
our case. This may
happen in infinite of finite time depending upon initial data. If this event
occurs in finite time, then subsequently, the solution  shrinks to  a point.
This behavior can be illustrated by the pictures below. The precise
meaning is contained in Theorem \ref{sant1}.

\begin{picture}(100,140)(-280,-80)
%%%%%%%%%%%%%%%%%%%%%%%%%%%%%%%%

%%% RYSYNEK1
\put(-200,0){.}
%punkt(-10,50)
\qbezier(-210,50)(-235,45)(-250,10)
\qbezier(-190,50)(-200,53)(-210,50)

%punkt(-50,-10)
\qbezier(-250,-10)(-235,-45)(-210,-50)
\qbezier(-250,10)(-253,0)(-250,-10)

%punkt(-50,10)
\qbezier(-190,-50)(-170,-46)(-150,-10)
\qbezier(-210,-50)(-200,-53)(-190,-50)

\qbezier(-150,-10)(-147,0)(-150,10)

%punkt(10,50)
\qbezier(-190,50)(-183,50)(-185,45)

\qbezier(-185,45)(-188,33)(-185,20)

%punkt(50,10)
\qbezier(-150,10)(-152,17)(-155,15)

\qbezier(-155,15)(-167,12)(-180,15)

%puntk(20,15)
\qbezier(-180,15)(-185,15)(-185,20)

\put(-220,-70){Time $t=0$}

%%%%%%%%%%%%%%RYSUNEK 2
%%%%%%%%%%%%%
\put(0,0){.}
%punkt(-10,50)
\qbezier(-10,50)(-15,45)(-50,10)
\qbezier(-10,50)(0,53)(10,50)

\put(-28,28){\vector(1,-1){3}}

%punkt(-50,-10)
\qbezier(-50,-10)(-45,-15)(-10,-50)
\qbezier(-50,10)(-53,0)(-50,-10)

\put(-28,-28){\vector(1,1){3}}

%punkt(-50,10)
\qbezier(10,-50)(15,-45)(50,-10)
\qbezier(-10,-50)(0,-53)(10,-50)

\put(28,-28){\vector(-1,1){3}}

\qbezier(50,-10)(53,0)(50,10)

%punkt(10,50)
\qbezier(10,50)(11,49)(15,45)

\qbezier(15,45)(12,33)(15,20)

\put(10,46){\vector(-1,-1){3}}

%punkt(50,10)
\qbezier(50,10)(49,11)(45,15)

\qbezier(45,15)(33,12)(20,15)

\put(46,10){\vector(-1,-1){3}}

%puntk(20,15)
\qbezier(20,15)(19,16)(15,20)

\put(19,19){\vector(1,1){3}}

\put(-20,-70){Time $t=t_1$}
\end{picture}

%%%%%%%%%%%%%%%%%%%%%%%%%%%%%%%%%%%%%%%%%%%%%%%%%%%%%%%%%%%

The evolution is determined by motion of facets defined by
singularities of the $J$-function (the arrows show the direction of the
evolution). In finite time we obtain a convex domain,
which becomes a square converging to a point in finite time.

%%%%%%%%%%%%%%%%%%%%%%%%

\begin{picture}(100,140)(-100,-80)
%%%%%%%%%%%%%%%%%%%%%%%%%%%%%%%%
 %%%%%%%%%%%%%%RYSUNEK 3%%%%%%%%%%%%%

\put(0,0){.}

\put(0,50){\line(-1,-1){50}}
\put(-50,0){\line(1,-1){50}}
\put(0,-50){\line(1,1){50}}

\put(0,50){\line(1,-1){14}}
\put(50,0){\line(-1,1){14}}

\put(15,25){\line(1,-1){10}}

\put(-24,24){\vector(1,-1){3}}
\put(-24,-24){\vector(1,1){3}}
\put(24,-24){\vector(-1,1){3}}
\qbezier(14,36)(12,33)(15,25)
\qbezier(25,15)(33,12)(36,14)
\put(6,42){\vector(-1,-1){3}}

%punkt(50,10)

\put(42,5){\vector(-1,-1){3}}

%puntk(20,15)

\put(22,22){\vector(1,1){3}}

\put(-25,-70){Time $t=t_2$}

%%%%%%%%%%%%%%%%%RYSUNEK 4

\put(150,0){.}

\put(150,30){\line(-1,-1){30}}
\put(150,30){\line(1,-1){30}}
\put(150,-30){\line(1,1){30}}
\put(150,-30){\line(-1,1){30}}

\put(164,14){\vector(-1,-1){3}}
\put(164,-14){\vector(-1,1){3}}

\put(136,14){\vector(1,-1){3}}
\put(136,-14){\vector(1,1){3}}

\put(125,-70){Time $t=t_3$}

%%%%%%%%%%%%%%
%% RYS 5
\put(250,0){.}
\put(225,-70){Time $t=t_{end}$}
\end{picture}

%%%%%%%%%%%%%%%%%%%%%%%%%%%%%%%%%%%%%%%%%%%%%%%

\noindent

We have to underline that the illustrated evolution hides the novel idea
of definition of singular term $\delta_0(u_x)u_{xx}$ being a multiplication of
two Dirac deltas (as in (\ref{tvf})),
however the nonlocal character   will allow us to define this
object. Additionally,  by the uniqueness of solutions to our system we show that
our novel definition is the only admissible. Formally, the dissipation
caused the Dirac delta coefficient is so strong that the changes of regularity
(i.e. appearance of the facets) happen instantly. % finite time.
%Of course  we obtain  ``singular'' end of facets, but maybe the purely
%classical point of view for this type of system is not the most appropriate.
%\smallskip [Moze to jest za skormnie :)]
%\smallskip

We will state our results in the Section below, the proofs will be presented in
the further Sections. Here, we present the outline of the rest of the paper.
We  show the existence of weak solutions   in Section \ref{sdwa},
uniqueness is the content of that Section, too. The
qualitative analysis is based on
the semi-discretization which is performed in Section \ref{sla}. Our
goal is to make some of the properties more apparent.  Namely, we want to show
that
facets (i.e. intervals where $\vfi=\La_s+s$ has a constant value equal to one
of the $\al_i$'s) form instantaneously.
In Section \ref{sas}, we show further geometric properties of solutions,
namely the curve becomes convex
(i.e. the angle $\vfi$ becomes monotone) in finite time. In addition, we show
that solutions become fully faceted in finite time, i.e. the solution is
composed only of facets.  These two events are
not correlated in time. Finally, we show that our solutions converge to a
special solution which we call minimal.

%%% --- s 2 ---
\section{The main results}\label{sdwa}

Here, we  present our results. We begin by noticing that
if $J$ is given by  (\ref{jot}), then the meaning of (\ref{rn2}) is not clear
at all because its right-hand-side formally becomes
$$
\La_t =
\frac\pi2\sum_{k=0}^3 \de_{k \frac{\pi}2 - \frac{3\pi}4}(s+\La_s)
\La_{ss}.
$$
Hence the above equation can be viewed as  a generalization of equation (\ref{tvf}).

We will use the tools of the convex analysis to interpret it. Due to
convexity of $J$ its subdifferential is always well-defined. Since in general
$\partial_\phi J(\phi)$ is not a singleton it is necessary to find its proper
selection, in particular (\ref{rn2}) takes the form,
\begin{equation}\label{nasz}
\begin{array}{ll}
\displaystyle \La_t \in
\frac{\partial}{\partial s}
{\partial_\vfi J(\La_s+s),} & \hbox{in }S\times(0,T),\\
\La(s,0) =\La_0(s), &\hbox{on } S,\\
 \La(2\pi,t) = \La(0,t), & \hbox{for } t\geq 0,
\end{array}
\end{equation}
where $S$ is the unit circle.

In other words, we have to find (weakly) differentiable selections of
$\partial_\phi J(\La_s +s)$.  Thus, we are lead to the following
notion of a weak solution to (\ref{rn2}).

\begin{definition}\label{defi}
{\it We say that  $\La\in
C([0,T];L_2(S %0,2\pi
))$,
such that $\La_s\in  L_\infty(0,T;\TV(S %0,2\pi
))$ is {\it a weak solution} to (\ref{nasz}), if there exists a function
$\Om\in L_1(0,T;W^{1}_{1}(S))$ such that
$\Om(s,t)\in \partial I(\La_s+s)$ a.e., and %$\Om_s(\cdot,t)$
%is integrable over $S$ %[0,2\pi)$ such that
for any function $h$ %smooth, $2\pi$ periodic function $h$
in  $C^\infty(S)$
it holds
$$
%\int_{\al_0}^{\al_0+2\pi} \La_t h = -\int_{\al_0}^{\al_0+2\pi} \Om
%h_s.
\int_S \La_t h = -\int_S (\Om-s)  h_s +\int_S h.
$$}
\end{definition}

With this definition we can show the following existence result.
\begin{theorem}\label{ist}
{\sl  Let us suppose that $J$ is defined by(\ref{jot}),
     $\La_0 \in L_1(S)$ and $\La_{0,s}\in \TV(S)$ , then
    there exists $\La\in
    C^\al(0,T;L_2(S))$ with $\alpha>0$, additionally
    $$
    \La_s\in
    L_\infty(0,T;\TV(S)) \mbox{ \ \  and  \ \ } \La_t\in
    L_2(0,T;L_2(S%0,2\pi
))
    $$
     such that it is a unique weak solution to  (\ref{nasz}).
}
\end{theorem}

The proof will be achieved through an approximation procedure, it is performed
in Section 3. Moreover, we  show uniqueness of the solution constructed here,
this is the content of Theorem \ref{twjed} in Section 3.

However, our main goal is to describe precisely qualitative properties
of solutions to
(\ref{nasz}). As a motivation, we present a special type of solutions, which we
will call {\it minimal solutions}, which are given explicitly, one of them is
given here, (see also \S 3.2),
$$
\bLa (s,t) = \int_0^s\bfi(u)\,du + t,
$$
where
\begin{equation}\label{rn-min}
\bfi(s) = \frac\pi4 \chi_{[0, \frac\pi2)}(s) +
\frac{3\pi}4 \chi_{[\frac\pi2,\pi )}(s) +
\frac{5\pi}4 \chi_{[\pi,\frac{3\pi}2 )}(s) +
\frac{7\pi}4 \chi_{[\frac{3\pi}2,2\pi)}(s).
\end{equation}
It is a matter of an easy exercise to see that $\bar\La$ defined above with
$\bar\Om(x,t) =x$ is indeed a weak solution
to (\ref{nasz}). In fact, this an asymptotic profile, which can be
reached in finite time.
%Employing our geometrical interpretation (as on pictures in Introduction)
%the minimal solution is a squares [kurczacy sie] to a point in finite time.

We will keep in mind this example while developing the proper class of regular
solution. %Main quantity is the obtained regularity of solutions and meaning
	  %in which  the solutions fulfill the system.
The idea is that we want to extend properties of
convex solutions to a more general class, hence we introduce a class of
J-regular function, where restrictions on regularity depend  on function $J$
from (\ref{jot}).

Firstly, we define the space of functions which are helpful to describe the
regularity of the derivative of our solutions. We recall that
any function $\phi \in TV$ is a difference of two monotone
functions. Thus, we shall call a multifunction
$\phi:[0,2\pi) \to 2^{\bR}$ a {\it maximal $TV$ function} if it is a
difference of two maximal monotone multifunctions and one of them is
continuous.

\smallskip

\begin{definition}\label{def2} {\sl
We say that a maximal $TV$  multivalued
function $\phi :[0,2\pi) \to \bR$ is J-regular,
i.e. $\phi \in \JR [0,2\pi)$, provided that
%iff  each selection of  $\phi \in TV[0,2\pi)$ and
the set
$$
\Xi(\phi)= \{ s \in [0,2\pi): \phi(s) \ni \alpha_k
  \mbox{ \ for \ } k=0,1,2,3\}
$$
consists of a \underline{finite number} of  connected components, i.e. we allow
only isolated  intervals or isolated points. Additionally, on any
connected subset
$[0,2\pi)\setminus \Xi$ function $\phi$ takes its values in interval
$(\alpha_k, \alpha_k+\frac{\pi}{2})$ for
  some $k=0,$ \dots, 3, modulo $2\pi$ -- see (\ref{angle}).

For each $\phi \in \JR [0,2\pi)$ we define a function
$K:\JR [0,2\pi) \to \mathbb N$ by the formula
$$
K(\phi)= %\{
\mbox{the number of connected components of the set $\Xi(\phi)$}.
%\}
$$
Additionally we put
$$
||\phi||_{\JR [0,2\pi)}=||\phi||_{TV[0,2\pi)}+K(\phi).
$$}
\end{definition}
\smallskip

Let us note that the $\JR$ class does not form a Banach space. It is
not a linear space. In order to formulate the meaning of solutions,
first we define the {\it composition  of $\JR$ functions with
  $\partial J$}. Because of the complex  structure the definition is long.

\begin{definition}\label{def3}
We define the {\it composition $\partial J\bar \circ A$,}
$$
\partial J \bar \circ A:[a,b]\to [e,f],
$$
where  $A:[a,b]\to [c,d]$ is an  $\JR$
function and
$\partial J:[c,d]\to [e,f]$
as follows:
\end{definition}

\smallskip

To begin with, we decompose  the domain $[a,b]$
 into three  disjoint parts
$[a,b]={\cal D}_r\cup {\cal D}_f \cup {\cal D}_s$,
where
\begin{equation}\label{de1}
\begin{array}{c}
{\cal D}_s=\{ s \in [a,b] : %A(s) \mbox{ \ is multivalued i.e. \ }
A(s)=[c_s,d_s] \mbox{ and } c_s <%\neq
d_s\}; \\[8pt]
{\cal D}_f=\{ \bigcup_{k} (a_k,b_k): A|_{(a_k,b_k)}=c_k, \mbox{ where }
c_k \mbox{ is a constant}\};
%\mbox{ \ \ and }
\quad
{\cal D}_r=[a,b]\setminus ({\cal D}_s \cup {\cal D}_f).
\end{array}
\end{equation}
Then, the composition is defined in three steps:

\smallskip

1. For each $s\in {\cal D}_r$ the set  $A(s)$ is a  singleton, thus
the composition is given in  the classical way
\begin{equation}\label{de2}
\partial J\bar \circ A(s)=\partial J(A(s)) \mbox{ \ \ \ \ \ \ for \ } s\in {\cal D}_r.
\end{equation}

2. In the case $s\in{\cal D}_f$ the definition is ``unnatural''.
For a given set $(a_k,b_k)\subset {\cal D}_f$ we have
$A|_{(a_k,b_k)}=c_k$.
%\smallskip
%
If $\partial J(c_k)$ is single-valued, then for $s\in (a_k,b_k)$ we have,
$$
\partial J\bar \circ A(s)=\{\frac{dJ}{d\phi}(c_k)\}.
$$

\smallskip

However, if $\partial J(c_k)$ is multivalued, i.e. $\partial
J(c_k)=[\alpha_k,\beta_k]$, then the definition is not immediate.
We have to consider four cases related to the behavior of multifunction $A$ in a
neighborhood of interval $(a_k,b_k)$. The regularity properties  of the
$\JR$ class imply the necessity to consider the
following four cases  (for small $\epsilon>0$):

\smallskip

(i) \ $A$ is increasing, i.e. $A(s) < c_k$ \ for $ s \in (a_k-\epsilon,a_k)$
and $A(s)> c_k$ \ for $s \in (b_k,b_k+\epsilon)$;

\smallskip

(ii) \  $A$ is decreasing, i.e. $A(s) > c_k$ \ for $ s \in (a_k-\epsilon,a_k)$
and $A(s)< c_k$ \ for $s \in (b_k,b_k+\epsilon)$;

\smallskip

(iii) $A$ is convex, i.e. $A(s) > c_k$ \ for $ s \in (a_k-\epsilon,a_k)$
and $A(s)> c_k$ \ for $s \in (b_k,b_k+\epsilon)$;

\smallskip

(iv) $A$ is concave, i.e. $A(s) < c_k$ \ for $ s \in (a_k-\epsilon,a_k)$
and $A(s)< c_k$ \ for $s \in (b_k,b_k+\epsilon)$.

\smallskip

\noindent
The case (i) we put
\begin{equation}\label{de3}
\partial J \bar \circ A(t)=x_k(t-b_k)+y_k(t-a_k)
\mbox{ \ \ \ \ \ \ for \ } t\in (a_k,b_k),
\end{equation}
where
$x_k=\frac{\alpha_k}{a_k-b_k}$  and  $y_k =\frac{\beta_k}{b_k-a_k}$.

For case (ii) we put
\begin{equation}\label{de4}
\partial J\bar \circ A(t)=x_k(t-b_k)+y_k(t-a_k)
\mbox{ \ \ \ \ \ \ for \ } t\in (a_k,b_k),
\end{equation}
where
$x_k=\frac{\beta_k}{a_k-b_k}$  and  $y_k =\frac{\alpha_k}{b_k-a_k}$.

When we deal with case (iii) we set
\begin{equation}\label{de5}
 \partial J \bar \circ A(t)=\beta_k
\mbox{ \ \ \ \ \ \ for \ } t\in (a_k,b_k).
\end{equation}

Finally, if  (iv) holds, then we put
\begin{equation}\label{de6}
\partial J \bar \circ A(t)=\alpha_k
\mbox{ \ \ \ \ \ \ for \ } t\in (a_k,b_k).
\end{equation}

3. In  the last case, if $s\in{\cal D}_s$ our definition is just a
consequence of  first two steps. Since  set ${\cal D}_s$ consists of a
countable number of points we consider each of them separately. We have
$A(d_k)=[e_k,f_k]$ with $e_k\neq f_k$,
then
\begin{equation}\label{de7}
 \partial J \bar \circ A(d_k)=[\limsup_{t \to d_{k}^-}  \partial J \bar \circ A(t),
\liminf_{t\to d_k^+} \partial J \bar \circ A(t)].
\end{equation}
Definition \ref{def3} is complete.

\smallskip

Thanks to the $\JR$ regularity of $A$, the above limits are well
defined. As a result,
we are able to omit point from ${\cal D}_s$ in (\ref{de1}). We note
that the above
construction guarantees that
$$
\partial J \bar \circ A:[a,b]\to [e,f] \mbox{ \ \ is a $\JR$
function.}
$$

\smallskip

After having completed the definition we make
additional comments on step 2. Formulae (\ref{de3})-(\ref{de6}) are
immediate consequences of
the pointwise %follow straightforward from
approximation of the considered function by smooth functions.
The presented composition  agrees with the results from
\cite{MuchaRybka}, where a stationary version of the
problem has been considered. In particular, our definition follows from
a requirement: if $A$ is maximal monotone then we expect
$$
A^{-1}\bar\circ A=Id.
$$
Moreover, the composition of two maximal increasing functions is
maximal increasing. Another point, which should be emphasized, is  the
nonlocal character of the above
definition. Step 3 depends on step 2, so steps 1 and 2 should be
performed at the very beginning.

%\paragraph{Example.} It is easy to check that $\bar\Om$, a part of the
%minimal solution, can be written as $\partial J\bar\circ \bar\phi$.
\bigskip

Now we are prepared to introduce the main definition.

\smallskip
\begin{definition}\label{def4}
{\it We say that a function $\Lambda: S \to \bR$ is  an
\underline{almost classical solution} to system (\ref{rn2}) iff
$\La$ is a weak solution with $\Om=\partial J\bar \circ
          [\Lambda_{s}+s]$, $\Lambda_{s}+s \in  L_\infty(0,T;\JR
          [0,2\pi))$ and
\begin{equation}\label{de9}
\begin{array}{lcr}
\Lambda_t=\frac{d}{ds} \partial J\bar \circ [\Lambda_{s}+s]&
\mbox{ in } & [S \times ((0,T)\setminus N)]\setminus \bigcup_{0< t< T} \partial
\Xi({\Lambda_{s}(\cdot,t)+s})\times\{t\},\\
\Lambda|_{t=0}=\Lambda_0 & \mbox{ on } & S,
\end{array}
\end{equation}
where $N$ is finite and $\partial E$ %\Xi({\Lambda_{,s}(\cdot,t)+s})$
denotes the boundary of set $E$.
%\Xi$ for function $\Lambda_{,s}(\cdot,t)+s$.
}
\end{definition}

%[!!!!!!!!! tu trrzbe sprawdzic co znaczy pochodna czasowa !!!!!!!]

\smallskip
The main point of
Definition \ref{def3} is to determine the composition appearing on the
RHS of the equation on sets,
where the solution and $\partial J$ are singular.
Note that equation $(\ref{de9})_1$ is fulfilled in the classical sense
except for finite number of point
for each $t \in (0,T)\setminus N$. This is so due to the definition of set
$\Xi({\Lambda_{s}(\cdot,t)+s})$ implying that its boundary
consists of finite number of points. It is easy to see that the
minimal solutions (\ref{rn-min}) fulfills Definition \ref{def4}.

The main result of our considerations is the following.

\begin{theorem}\label{thmreg} {\sl
Let $\Lambda_0$ be such that $\Lambda_{0,s}+s \in
\JR [0,2\pi)$,
then there exists a unique almost classical solution to system (\ref{nasz})
conforming to Definition \ref{def4}.}
\end{theorem}

In fact this is a statement about the regularity of weak solutions.
%The proof of
Theorem \ref{thmreg} is a result %will follow from the analysis
of the
semi-discretization of system (\ref{rn2}). At this level, we will
be able to show that facets  must appear, as suggested by the pictures
in the Introduction. The
semi-discretization will determine  the RHS of (\ref{de9}) on
sets where the solution falls into the singular part
of $\partial J$. We will obtain that on these sets the term $\partial
J$ is constant
on each connected part (or time dependent for the evolutionary system). Next,
by the elementary means we will show that the semi-discretization
tends uniformly to the
solutions obtained by Theorem \ref{ist}.  However, performing a rigorous proof
that we indeed constructed an almost normal solution requires more work on the
structure of weak solutions, which is the content of Section
\ref{sas}. Thus, it will be postponed until the end of this part.

At the end, in Section \ref{sas}, we deeply go into the
qualitative analysis of the evolution showing the convexification effect and
convergence to the minimal solutions. Since we know that facets must
appear and the solutions are unique we are in a position  to construct
quite explicit solutions. We are able to follow their qualitative
changes. This is made precise in Theorem \ref{sant1}.
In particular we show instantaneous creation of facets. For the sake of
this study we show a comparison principle in subsection
\ref{sas}.1. Moreover, we show that the evolution of facets is governed by a
system of ODE's which  are coupled if the facets interact, this is
explained in Section \ref{sas}. A conclusion from our analysis is
existence of a sequence
of instances at which our solution gets simplified before it gets the
final form of the asymptotic profile, i.e. the minimal solution.

%%% --- s 3 ---
\section{Existence of solutions}\label{sla}

In this Section we  show an existence and uniqueness of weak solutions of
%result for
(\ref{rn2}). We  use the tools of the convex analysis to interpret it.  In
particular,
we shall make the gradient flow structure of  (\ref{rn2})
transparent.  However, the
existence is shown by the method of regularization.
Some of the statements are easier to interpret if they are written in the
language of the `angle' $\vfi= \La_s +s$. Here, $\vfi$ plays the role of the
angle between the normal to the curve and the $x_1$--axis. Thus, for convex
closed curves $\vfi$ must be increasing, but we shall not require
that, instead we admit $\vfi$ being a functions of bounded total
variations, i.e., $\vfi(\cdot,t)\in TV(S)$, in particular $\vfi\in
L^\infty(S)$ and it may be discontinuous though.

%%% --- s 3.1 ---

\subsection{The proof of the general existence result}

%While the example of the minimal solutions is not very much instructive
%as far as the time regularity is concerned, nonetheless it shows that
%(\ref{genw}) is satisfied. It also gives us a motivation for our notion  of
%solution given in Section 2. Here we notice an advantage of this
%definition. Namely,  the first integrand on the RHS has zero average
%over $S$.

%Let us notice that for any $a\in(-\frac34\pi,-\frac\pi4)$ we have $\partial
%J(a)=\{-\frac\pi2\}$, hence by the definition of $J$ it follows that
%$J(a+2\pi)=\{-\frac\pi2+2\pi\}$. As a result, the integration by parts, which
%is permitted for the class of admissible function yields $\int_S \frac
%  {d\Om}{ds}$ on the RHS.

We present a proof of our existence result, Theorem \ref{ist}. It will be
achieved through an approximation procedure.
For any $\ep>0$ we set
\begin{equation}\label{jotep}
J^\ep(x):= J\star \rho_\ep(x)+\frac{\ep^2}2 x^2,
\end{equation}
where $\rho_\ep$ is a standard mollifier kernel, with support in
$(-\ep,\ep)$. Let us note properties of the approximation $J^\ep$:

(a) $J^\ep\in C^\infty(\bR)$;\\

(b) $\frac d{dx} J^\ep$ is strictly monotone; \\

(c) $\frac {d^2}{dx^2} J^\ep\ge\ep$;\\

(d) $\frac d{dx} J^\ep(x) -\ep x = \frac d{dx} J(x)$ for
$x$ such that $|x-\al_k|>\ep$ for
$k=0,1,2,3$.

\medskip\noindent
We start with existence of the regularized system.
\begin{lemma}\label{lem21}{\sl Let us suppose that $J^\ep$ is defined by
(\ref{jotep}) and $\La_0^\ep$ is smooth and $2\pi$-periodic.
Then, for any $T>0$ there exists a unique, smooth solution to
the regularized problem,}
\begin{eqnarray}\label{naszr}
&&\La^\ep_t = %\frac d{ds}
\frac{\partial}{\partial s}
\frac d{d\vfi}J^\ep(\La^\ep_s+s), \qquad\hbox{in }
S \times (0,T),\nonumber \\
&&\La^\ep(s,0) =\La_0^\ep(s),  \qquad\hbox{on } S,\\
&& \La^\ep(s+2\pi,t) = \La^\ep(s,t),  \qquad\hbox{for }\quad
 t>0.\nonumber
\end{eqnarray}
\end{lemma}

\medskip\noindent{\it Proof.}
By properties (a), (b) (c) and (d) of $J^\epsilon$, see (\ref{jotep}),
the existence and uniqueness of smooth solutions to (\ref{naszr}),
 is guaranteed by the standard theory of parabolic systems,
 see \cite{lady}. \qed

We now study properties of established solutions.

\begin{lemma}\label{lem32}
{\sl Let us suppose that $\La^\ep$ is a smooth  solution to
(\ref{naszr}).\\% and $\La^\ep$ is $2\pi$-periodic. \\
(a) If for $a,b \in {\mathbb R}$
and the initial datum satisfies $a\le (\La^\ep_{0,s}(s)+s)\le b$,
then, for all  $t<T$ we have
$$
a\le (\La^\ep_{s}(s,t)+s)\le b. \\
$$
(b) If moreover,
$(\La^\ep_{0,s}(s)+s)_s \in L_1(0,2\pi)$,
then, for all  $t<T$ we have
$$
 (\La^\ep_{s}(s,t)+s)_s \in L_\infty(0,T;L_1(0,2\pi)).
$$}
\end{lemma}

\medskip\noindent{\it Proof.} We use the maximum principle. First of
all, we differentiate (\ref{naszr}) with respect to $s$,
$$
\La^\ep_{st} = \frac d{ds}\left(
\frac{\partial^2
  J^\ep}{\partial\vfi^2}(s+\La^\ep_{s})(s+\La^\ep_{s})_s
\right).
$$
We notice $(s+\La^\ep_{s})_t = \La^\ep_{st} $. We set
$w =(s+\La^\ep_{s})$, hence we obtain the equation for $w$,
\begin{equation}\label{m1}
w_t = \frac d{ds}(a(s,t)w_s),
\end{equation}
where by (\ref{jotep})
we have
$a(s,t) =
\frac{\partial^2 J^\ep}{\partial\vfi^2}(s+\La^\ep_{s})\ge\ep>0$.
Hence, by the maximum principle we obtain $(a)$.

To prove (b) we note that from (\ref{m1}) we obtain
\begin{equation}\label{m2}
w_{st}=\frac{d^2}{ds^2}(a(s,t)w_s).
\end{equation}
By Lemma \ref{lem21} our solutions are  smooth. In oder to finish the proof of
(b)  it is enough to integrate (\ref{m2}) over
sets $\{ w_s >0\}$ and $\{ w_s <0\}$ to reach,
\begin{equation}\label{m3}
\frac{d}{dt}\int_{\{ w_{s}>0\} }w_{s} dx \leq 0
\mbox{ \ \ \ \ and \ \ \ \ } \frac{d}{dt}\int_{\{ w_{s}<0\} }w_{s} dx \geq 0.
\end{equation}\qed

Having established  this Lemma, we will obtain $L_\infty$ estimates for
the spatial derivative of solution $\La$.

\begin{corollary}\label{cor31}
{\sl There is a constant $M$ independent of $\ep$ and $T$ such  that}
$$
\|\vfi^\ep\|_{L_\infty(S\times(0,T))} \le M,\qquad
\|\vfi^\ep(\cdot,t)\|_{L_\infty(0,T;TV[0,2\pi))} \le M.
$$
\end{corollary}
\noindent{\it Proof.} The first part follows from  Lemma \ref{lem32} (a)
directly, because $\vfi^\ep = \La_s^\ep+s$. The
second part is the result  of Lemma \ref{lem32} (b), combined with the
properties of approximation of $TV$ functions in
$L_1$. %-space.  
\qed

We want to show that the estimates  for $\La^\ep$ will persist
after passing to the limit with $\ep$.

\begin{lemma}\label{lem33}
{\sl Let us suppose that $\La^\ep$ converges weakly in
$L^2(S %(0,2\pi)
\times(0,T))$ to $\La$. If $(\La^\ep_{s}+s)_s\ge0$ in
$\cD'(S %0,2\pi
)$, then $(\La_s+s)_s\ge 0$ as well in $\cD'(S)$.}
\end{lemma}
\noindent{\it Proof.} Indeed, if $h\in\cD(S)$ is positive, then
$0\le  \int_S (\La^\ep_s +s) h_s$.
The inequality holds after taking the limit.\qed %%znacznik  23.11.06

\begin{lemma}\label{lem34}
{\sl There is a constant independent of $\ep$ such that}
$$
\int_0^T\int_0^{2\pi} (\La^\ep)^2\,dxdt \le M,\qquad
\int_0^T\int_0^{2\pi} [(\La^\ep_x)^2+(\La^\ep_t)^2]\,dxdt \le M.
$$
\end{lemma}

\noindent{\it Proof.} The bound on $\int_0^T\int_0^{2\pi} (\La^\ep)^2$ is
trivial, due to $L^\infty$ estimates established in previous lemmas. Similarly,
the bounds in Corollary \ref{cor31} imply that
$\int_0^T\int_0^{2\pi} (\La^\ep_x)^2 \le M$.
We shall
calculate the last integral with the help of integration by parts,
\begin{eqnarray*}
\int_0^T\int_0^{2\pi}(\La^\ep_t)^2\,dsdt&& =
-\int_0^T\int_0^{2\pi}\La^\ep_{st}  \frac{d}{d\varphi}
J^\ep(\La^\ep_s+s)\,dsdt
+ \int_0^T \La^\ep_t \frac d{d\vfi}J^\ep(\La^\ep_s+s)|_{s=0}^{s=2\pi}\,dt\\
&&= \int_0^{2\pi} J^\ep(\vfi_0(s))\,ds - \int_S J^\ep(\vfi(s,T))\,ds\\
&&\quad
+ \int_0^T \La^\ep_t(0,t)\left( \frac d{d\vfi}J^\ep(\La^\ep_s(0,t)+2\pi) -
\frac d{d\vfi}J^\ep(\La^\ep_s(0,t)) \right)dt, 
\end{eqnarray*}
Here, we also exploited periodicity of $\La$. We notice that the
difference $\frac d{d\vfi}J^\ep(\La^\ep_s(0,t)+2\pi) - \frac
d{d\vfi}J^\ep(\La^\ep_s(0,t))$ equals exactly $2\pi$. Hence,
$$
\int_0^T\int_0^{2\pi}(\La^\ep_t)^2\,dsdt \le \int_0^{2\pi}
J^\ep(\vfi_0(s))\,ds +2\pi (\La^\ep(0,T) -\La^\ep(0,0))\le M
$$
due to the Corollary \ref{cor31}.  \qed

\medskip\noindent{\bf Remark.} 
We want to stress that the above estimate on $\La_t$ is one of the
most important differences  between (\ref{rn0}) and (\ref{nasz}). 

Now, we have enough information to select a weakly convergent
subsequence, with properties announced in the theorem.

%{\it Proof of Theorem \ref{ist}.} {\it Step 1.}
\begin{proposition}\label{pr31} {\sl There exists a
  subsequence $\{\ep_k\}$ converging  to zero, such that

(a) $\La^{\ep_k} \rightharpoonup \La$ in
  $W^{1}_{2}(S\times(0,T))$;\qquad
$\vfi^{\ep_k}_s \rightharpoonup \vfi_s $ as measures in $S\times (0,T)$.

(b) $\La \in C([0,T),L_2(S))$.}
\end{proposition}
\noindent{\it Proof.}  The first part of (a) is implied by
Lemma \ref{lem34}. The second part  of (a) follows from $\vfi^\ep =
\La^\ep_s +s$,  and Lemmas \ref{lem32}, \ref{lem34}.
Part (b) follows from Lemma \ref{lem32} and \ref{lem34}
 and the embedding theorem (we have already proved $\Lambda^\epsilon \in
 L_2(0,T;W^1_2(0,2\pi))\cap W^1_2(0,T;L_1(0,2\pi))$).
\qed

The next step is to show that  that the limit is indeed a solution. In
particular, we have to pass to the limit in the non-linear term.
First of all, we shall change the notation in order to make 
more transparent what we are doing. We want to find $w(s,t)$ such that
$w_s(s,t) =\vfi (s,t)$. By a simple integration of this formula and the
definition of $\vfi$, we can see
$$
w(s,t)= \frac12 s^2 + \La(s,t),
$$
where we set $w(0,t) = \La(0,t)$. Hence, $w_s = \vfi$ and we can
re-write the evolution problem as a gradient system
\begin{eqnarray}\label{rew}
&& w_t \in \frac{d}{ds}\partial J(w_s),\qquad\hbox{in }
S\times(0,T), \nonumber\\
&& w(s,0) = \frac12 s^2 + \La_0(s),\qquad  \mbox{ for } s\in S,\\
&& w(s,t) -\frac 12 s^2 \mbox{ \ is periodic for }\qquad t\in(0,T).\nonumber
\end{eqnarray}
If $\vfi(\cdot,0)$ is increasing, then due to Lemma
\ref{lem32} (b) and  Lemma \ref{lem33} $\vfi(\cdot,t)$ is increasing as well,
hence $w(\cdot,t)$ is
convex.
Obvious changes are required to write the system for the
regularization $w^\ep(s,t) =\frac12 s^2 + \La^\ep(s,t)$.

\begin{proposition}\label{pr32}
{\sl  For any fixed $t\ge 0$ and a sequence $\{\ep_k\}$ converging  to zero
there exists its subsequence
$\{\ep_k\}$ (not relabeled), such that for each $x\in[0,2\pi)$
the limit
$$
\lim_{\ep\to 0}\frac{d}{d \vfi}(J^\ep)(\vfi^\ep)(x,t) = \Om(x,t)
$$
exists. Moreover, $\Om(x,t)\in \partial J(\vfi(x,t))$ for almost every
$x\in[0,2\pi)$.}
\end{proposition}
{\bf Remark.} It is important for us to make the selection of the
subsequence independently of $t$.

\bigskip\noindent{\it Proof.} Indeed, once we fix $t>0$, we may recall that
$\vfi^\ep(\cdot, t)\in TV$ as well as $\frac{d}{d\vfi}J^\ep(\vfi^\ep)(x,t)\in
TV$. Hence, by Helly's convergence theorem there
exists a subsequence $\ep_k$ such that these sequences converge. Using the new
notation, we write,
$$
\lim_{\ep \to 0}\vfi^\ep(x,t)= w_x(x,t),\qquad
\lim_{\ep \to 0}\frac{d}{d\vfi}J^\ep(w^\ep_x(x,t)) = \Om(x,t).
$$

Now, we shall show that for each point $x$
the number $\Om(x,t)$ belongs to $\partial J(w_x(x,t))$. Since the
functions $J^\ep$ are convex, we have the inequality
$$
\int_0^{2\pi} J^\ep(w^\ep_x(x,t)+ h_x(x)) -J^\ep(w^\ep_x(x,t))\,dx
\ge\int_0^{2\pi} \frac{d}{d\vfi} J^\ep (w^\ep_x(x,t))  h_x(x)\,dx,
$$
for each $h \in C^\infty_0(0,2\pi)$.
We know that $w^\ep$ and $\frac{d}{d\vfi}J^\ep(w^\ep_x(x,t))$ have pointwise
limits, which are bounded, hence after passing to limit our claim will follow,
$$
\int_0^{2\pi} J(w_x(x,t)+ h_x(x))-J(w_x(x,t))\,dx
\ge\int_0^{2\pi} \Om(x,t)h_x(x)\,dx.\eqno\qed
$$
%[TEORETYCZNY KLOPOT, CIAG $\ep_k$ ZALEZY OD $t$, TYLE ZE GRANICA ISTNIEJE.]

%[NIESTETY, TRACIMY INFORMACJE PUNKTOWA. NIE WIEMY JESZCZE JAK PRZEJSC
%  DO GRANICY W WYRAZIE NIELINIOWYM].
We  finish the {\it proof of Theorem \ref{ist}.} By previous
Lemmas there exists a sequence $\La^\ep$ which converges weakly in
$W^{1}_{2}(S\times(0,T))$. In particular, if $h\in C^\infty_0(0,2\pi)$, $t>0$
and $\tau>0$ is arbitrary, then we  see
$$
\int_{t-\tau}^{t+\tau}\int_S \La^\ep_t h\,dsdt' =
\int_{t-\tau}^{t+\tau}\int_S
\frac{\partial }{\partial s}\frac{\partial }{\partial \vfi}
J^\ep(\La^\ep_s+s)h\,dsdt'
= -\int_{t-\tau}^{t+\tau}\int_S \frac{\partial }{\partial \vfi}
J^\ep(\La^\ep_s+s)h_s\,dsdt'.
$$
Since $\frac\partial{\partial\vfi}J^\ep(\La^\ep_s+s)$ is bounded, it
converges weak-$\ast$ in
$L_\infty((0,2\pi)\times(0,T))$ to $\Om$. We have to show that
$\Om(s,t)\in \partial J(\La_s+s)$. First we notice that we may
pass to the limit in the above integral identity,
$$
\int_{t-\tau}^{t+\tau}\int_S \La_t(s,t') h(s)\,dsdt'
=-\int_{t-\tau}^{t+\tau}\int_S \Om(s,t') h_s(s)\,dsdt'.
$$
By the Lebesgue differentiation theorem we deduce,
\begin{equation}\label{good}
\int_S \La_t(s,t) h(s)\,ds = - \int_S \Om(s,t) h_s(s)\,ds
\end{equation}
for a.e. $t\in [0,T]$ for $h \in W^1_2(S)$ (we used the fact that $0$ is not
distinguished on $S$).  In principle, the set
$G=\{t\in [0,T]:\ \hbox{ (\ref{good}) holds}\}$ depends upon $h$,
i.e. $G=G(h)$. We shall see, that in fact we can choose $G$
independently of $h$. Let us recall that $W^1_2(S)$ is separable
and let us suppose that $D$ is a dense, countable subset of
$W^1_2(S)$. Of course, ${\cal G} =\bigcap_{h\in D}^\infty G(h)$ is
a set of full measure. Let us then take $t\in{\cal G}$ and
$h\in C^\infty(S)$. Let us suppose that
$\{h_n\}$ is
a sequence in $C^\infty(S)$ converging to $h$ in the
$W^1_2(S)$-norm. Then,
$$
\int_S \La_t(s,t) h_n(s)\,ds = - \int_S \Om(s,t) (h_n)_s(s)\,ds
$$
for all $t\in{\cal G}$. We may pass to the limit with $n$ on both
sides, thus we reach,
$$
\int_S \La_t(s,t) h(s)\,ds = - \int_S \Om(s,t) h_s(s)\,ds.
$$
In other words, (\ref{good}) holds for all $h\in C^\infty(S)$ and all
$t\in{\cal G}$.

If we now fix $t\in{\cal G}$, we next apply Proposition \ref{pr32} to
deduce that
$\Om(s,t)\in\partial J(\La_s(s,t)+s)$. Hence the limit, $\La$, is indeed
a weak solution. \qed

Now, we are going to  prove uniqueness.

\begin{theorem}\label{twjed}{\sl If $\La^i$, $i=1,2$ are two solutions
with $\La^1(s,0)=\La^2(s,0)$, then $\La^1(s,t)= \La^2(s,t)$, for
$t\le  T$.}
\end{theorem}
\noindent{\it Proof.} If $\La^i$, $i=1,2$, are weak solutions, then by the
definition of weak solutions we have
$$
\int_S\La^i_t h\,ds =-\int_S (\Om^i-s) h_s\,ds+ \int_S h\,ds ,
$$
where $w^i \in-\partial J$ and $h$ is in $H^1$. We subtract these two
identities for $\La^2 $ and $\La^1$, then we take  $(\La^1- \La^2)$ as a the
test function. Finally, the integration over $(0,\bar t)$, $\bar t<T$
yields
$$
\int_0^{\bar t}\int_S\frac 12 \frac d{dt}(\La^1- \La^2)^2\,dsdt =
-\int_0^{\bar t}\int_S (\Om^1-\Om^2) (\La^1- \La^2)_s\,dsdt.
$$
Monotonicity of $\partial J$ implies that
$\frac12\|\La^1-\La^2\|^2_{L^2(S)} ({\bar t}) \le 0$.
Hence,
$\|\La^1-\La^2\|^2_{L^2(S)} ({\bar t})= 0$ for any $\bar t<T$.
\qed

%------- s 3.2 ----------------

\subsection{Minimal solutions}

It is well-known that important information about the studied system is
provided by special solutions, like traveling waves, self-similar
solutions and other symmetry solutions. We can not talk about
self-similar solutions because our systems  lacks direct geometrical
interpretation, however we may look for special ones, which we named
minimal solutions.

In the theory of curvature flows it is natural to anticipate existence
of curves such that their curvature is constant, but may change in
time. Here, we ask if there exists such a solution $\bfi$ to (\ref{rn2})
that
\begin{equation}\label{ink}
\frac d{d s}\partial J (\bfi)\ni k,  \mbox{ \ \ hence \ }
\partial J (\bfi) \ni ks +s^*,
\end{equation}
where $s^*$ is appropriately chosen, e.g. $s^*=\frac\pi4$ and $|k| =1
$. The last restriction is  of  geometric nature, namely we want that
for any $a\in\bR$ the image of $S$ %[a,a+2\pi)$
by $\partial J$ be contained in an interval no longer than $2\pi$.

In fact, we may come up with explicit formulas. One for $k=1$ is provided
by formula (\ref{rn-min}).
It is then obvious  that
\begin{equation}\label{roz}
\bfi(s) := (\partial J )^{-1}(s +s^*),
\end{equation}
as in \cite{MuchaRybka} and in Section 2. Moreover,
%\begin{equation}\label{rncal}
$\int_0^{2\pi}\bfi(s)\,ds = 2\pi^2 = \int_0^{2\pi} s\,ds.$
%\end{equation}
By the reversal of the orientation, we immediately obtain the solution for
$k=-1$,
$$
\bfi_{-1}(s) = \frac{7\pi}{4} \chi_{[0, \frac\pi2)}(s) +
\frac{5\pi}4\chi_{[\frac\pi2,\pi )}(s) +
\frac{3\pi}{4}\chi_{[\pi,\frac{3\pi}2 )}(s) +
 \frac\pi4 \chi_{[\frac{3\pi}2,2\pi)}(s).
$$
%
%Of course (\ref{rncal}) holds also for $\bfi_{-1}$.
%
We choose
$\bfi(s) := \bfi_{1}(s)$,  which is given by (\ref{rn-min}),
because we prefer to have $\bfi$ an increasing function.

As a result, $\bLa$ defined by
$\bLa (s,t) = \int_0^s\bfi(u)\,du + F(t)$
is indeed $2\pi$ periodic in $s$ and it is a solution to
(\ref{rn2}). Here, we must take $F(t) = A+t$. One
can check in  a straightforward manner that indeed $\bLa$ solves
(\ref{roz}). This is indeed so, because we have found
$\bfi(s)=\La_s(s)+s$ and $\Om$ is
a section of $\partial I(\bfi)$, namely, $\Om(s,t) =s$, which satisfies
(\ref{roz}). If we take $A=0$, then $\bar \La$ satisfies the initial
condition: $\bar\La(s,0) =\int_0^s\bar\vfi(u)\,du$.

%%%%%%%%%%%%%%%%%%%%%

%%%%--------s4----------

%%%%--------s5--------------

\section{The semi-discretization}% -- qualitative analysis}
\label{ssem}

In this part we examine the semi-discretization of (\ref{nasz}). Our goals are
not only to establish existence for the presented scheme, but also to show
qualitative properties of the obtained solutions. In
particular  our considerations will explain the appearance of facets. Finally,
we prove the convergence of 
solutions of the semi-discretization to the solutions obtained in Section 3.

We define the semi-discretization in time of system (\ref{nasz}) as follows
\begin{equation}\label{d1}
\frac{\lambda^k_h(s)-\lambda^{k-1}_h(s)}{h}\in \frac{d}{ds}\partial J[
\lambda^k_{h,s}(s)+s]
\end{equation}
and $\lambda^k_h(0)=\lambda^k_h(2\pi)$ and $(\lambda^0_h)_s=\phi_0$
for $k=1,\ldots,[T/h]$; or equivalently equation (\ref{d1}) can be stated
\begin{equation}\label{d7}
\lambda^k_h(s)-h\frac{d}{ds}\partial
J[\lambda^k_{h,s}(s)+s] \ni\lambda^{k-1}_h(s).
\end{equation}

We  establish existence of solution to this problem.

\begin{lemma}\label{provx}\quad {\sl
Let us suppose that an absolutely continuous function $v$ is such that  $v_s =
\vfi \in TV[0,2\pi)$, then there exists $u\in AC([0,2\pi))$ such that
$u_s\in TV$, which is a  solutions to (\ref{d1}), i.e. 
\begin{equation}\label{mm1}
u -v \in h\frac d{ds}\partial J(u_s)
\end{equation}
with $u(0)+\frac12 (2\pi)^2=u(2\pi)$ and the following bound is valid
\begin{equation}\label{mm1a}
||u_{s}||_{TV}\leq ||v_s||_{TV}.
\end{equation}
}
\end{lemma}

\bigskip\noindent{\bf Remark.} Our understanding of (\ref{mm1}) is the same
as that of (\ref{nasz}), i.e., there exists $\om \in W^1_1([0,2\pi))$, such
that $\om(x)\in \partial J(u_s)$ and $ u-v = h \frac{d}{ds}\om$.

We also note that $u$ and $v$ appearing in this Lemma need not be periodic,
on the other hand $\La(\cdot, t)$ and  $\lambda^k_h(\cdot)$ are periodic.

\bigskip\noindent{\it Proof.}
Let us notice that if $u$ is a solution to  (\ref{mm1}), then $0$ belongs to
the subdifferential of the functional
$$
\cJ(u) =
\int_0^{2\pi}[hJ(u_s) +\frac12(u-v)^2],\quad \hbox{for } u\in AC([0,2\pi), \
  u_s \in TV. 
$$
i.e. $u$ is a minimizer of $\cJ$. To be precise, we define $\cJ$ on $L_2(S)$
by the above formula for $u\in AC([0,2))\pi$ with $u_s\in TV$ and we put
  $\cJ(u)=+\infty$  for $u$ belonging to the complement of this set.

In order to solve (\ref{mm1}), we consider a family of regularized problems,
$$
\cJ_\ep(u) =
\int_0^{2\pi}[hJ_\ep(u_s) +\frac12(u-v)^2],
$$
where $J_\ep$ is the same regularization of $J$ that we used in (\ref{jotep}).

The functional $\cJ_\ep$ is well-defined, convex and coercive on
the standard Sobolev space $W^2_1(0,2\pi)$, thus it
possesses a unique minimizer $u^\ep$. Now, we apply again the methods used in
Section 3.1 to show existence of a weak solution of the evolution problem
(\ref{rn2}). The regularization of system (\ref{mm1}) leads to the following
equation
$$
u^\ep_{ss}-\frac{d^2}{ds^2}
(\frac{\partial J_\epsilon}{\partial \varphi^2}(u_s^\ep)u_{ss}^\ep)
=v_{ss}^\ep.
$$
By repeating the argument for (\ref{m2}),we get
$||u_{ss}^\ep||_{L_1}\leq ||v_{ss}^\ep||_{L_1}$.
Passing to the limit with $\epsilon \to 0$ yields (\ref{mm1a}).

In addition we have the following bounds
$\int_0^{2\pi} (u^\ep)^2\,dx \le M$,
$\int_0^{2\pi} (u_x^\ep)^2\,dx \le M$.
In order to prove them we follow the lines of reasoning of Corollary
\ref{cor31} and Lemma \ref{lem34}.
These bounds suffice to show existence of  a
  subsequence $\{\ep_k\}$ converging  to zero, such that\\
(a) $u^{\ep_k} \rightharpoonup u$ in
  $W^{1}_{2}(0,2\pi)$;\qquad
$u^{\ep_k}_{ss} \rightharpoonup u_{ss} $ as measures.

Subsequently, by Helly's theorem we conclude existence of the pointwise limits
(for another subsequence $\{\ep_k\}$, not relabeled) 
% such that for each $x\in[0,2\pi)$ the limit
$$
\lim_{t\to \infty}\vfi^\ep(x)= u_x(x),\qquad
\lim_{t\to \infty}\frac{d}{d\vfi}J^\ep(u^\ep_x(x)) = \Om(x).
$$
Moreover, $\Om(x)\in \partial J(u_x(x))$ for each $x\in[0,2\pi)$.

Now, we show uniqueness of solutions, constructed in Lemma \ref{provx}.

%-------lemat x.x--------
\begin{lemma}\label{lem4xx}\quad {\sl
Let $v\in AC([0,2\pi))$, $v_s\in TV([0,2\pi))$, then there exists at most one
weak solution $u\in AC([0,2\pi))$, $u_s\in TV(S)$ to problem (\ref{mm1})}
\end{lemma}
{\it Proof.} Let us suppose that there are  two
solutions to (\ref{mm1}), $u^i$, $i=1,2$. By the definition, there are
two functions
$\omega_i\in \partial J(u_s^i)$, $i=1,2$, such that
$$
u^i -v =h \frac d{ds}\omega_i,\qquad i=1,2.
$$
After subtracting these two equations and multiplying them by
$u_1-u_2$ and integrating over $[0,2\pi)$ we see
$$
\|u^1-u^2\|^2 - \int_0^{2\pi} h(\frac d{ds}\omega_1 -\frac d{ds}\omega_2)
(u^1-u^2)\,ds =0.
$$
The integration by parts leads us to
$$
0=\|u^1-u^2\|^2 +
\int_0^{2\pi} h(\omega_1 -\omega_2)(u^1_s-u^2_s)\,ds
\ge \|u^1-u^2\|^2 \ge 0 .
$$
As a result $u^1=u^2$.
\qed

In order to finish our preparations, %for the existence result 
we introduce the
sets of preferred orientation which %play 
dominate the behavior of
solutions. Let us suppose, that $w$ is absolutely continuous and $w_s\in TV$,
then at any point $s$, the left derivative $w_s^-$, as well as the right
derivative $w_s^+$  are well-defined, hence we may set
\begin{equation}\label{depaw}
\partial w(s) = \{\tau w_s^- + (1-\tau) w_s^+:\ \tau\in [0,1]\}.
\end{equation}
If $w$ is convex, then $\partial w$
is the well-known subdifferential of $w$.

Now, for each $l=0,1,2,3,$ we set
\begin{equation}\label{d4}
\Xi_l(w_s)={\{ s \in [0,2\pi]: w\hbox{ is differentiable at }s\hbox{ and }
w_{s}(s)=\alpha_k \quad\hbox{or}\quad
 \alpha_k \in \partial w(s)
%[w_{s}^-(s),w_{s}^+(s)]\quad\hbox{or}\quad
% \alpha_k \in [w_{s}^+(s),w_{s}^-(s)]
\}}
\end{equation}
Furthermore, we set
$\Xi(w_s)=\bigcup_{l=0}^3 \Xi_l(w_s)$.

The result, delivering the main properties of solutions, is the following.

\smallskip

%----------tw. x.x-------------
\begin{theorem}\label{thm4xx}\quad {\sl
Let $\phi_0=\lambda^0_{h,s}+s \in \JR [0,\pi)$.
Then a solution $\{\lambda^h_k\}$ to problem (\ref{d1}) exists, it is
unique and it satisfies the following
bound
\begin{equation}\label{d3}
||\lambda^h_{k,s}+s||_{\JR [0,2\pi)}\leq ||\lambda_{0,s}+s||_{\JR [0,2\pi)}.
\end{equation}
Moreover,
we have
\begin{equation}\label{d5}
\Xi(\lambda^h_{k-1,s}+ s) \subset \Xi(\lambda^h_{k,s}+ s)
\mbox{ \ \ and \ \ }
K(\lambda^k_{h,s}+s)\leq K(\lambda^{k-1}_{h,s}+s)
\end{equation}
and
\begin{equation}\label{d5a}
\sup_k \sup_{l=0,1,2,3}|\Xi^k_l\setminus\Xi^{k-1}_l|\leq C(V(h)+h^{1/2}).
\end{equation}
where $V(s) \to 0$ as $s \to 0$ and $V$ is determined by the initial datum
$\phi_0$. Moreover, 
on connected components of the set $\Xi^{k-1}\setminus (\bigcup_{l=0}^3
\Xi^k_l \setminus \Xi^{k-1}_l)$
\begin{equation}\label{d5b}
\frac{d}{ds} \partial J[\lambda^k_{h,s}+s] \mbox{ \ \ is constant.}
\end{equation}}
\end{theorem}

\smallskip\noindent{\it Proof.}
By Lemmas \ref{provx} and \ref{lem4xx} we conclude existence of the sequence
of solutions to the semi-discretization, the solutions are such that $\la^h_s$
belong to $TV(S)$.
It is enough to  restate the equation (\ref{d7})
 as follows:
\begin{equation}\label{q1}
u-h\frac{d}{ds}\partial J[u_{s}]=v
\end{equation}
with $u=\lambda^k_{h}+\frac 12 s^2$ and $v=\lambda^{k-1}_h+\frac 12 s^2$,
and boundary condition $u(0)+2\pi^2=u(2\pi)$.

The set, where function
$J[u_s]$ is singular, i.e.  $\Xi(u_s)$, plays the key role. 
%most important information will be described on a set where function
%$J$ is singular, i.e. on $\Xi(u)$.
Our first task is to prove the inclusion from (\ref{d5}).
Note that in a neighborhood of any point $s\notin \Xi(u_s)$ function
$\partial J[u_{s}(\cdot)]$ is constant, hence we get
$u(s)=v(s)$. Thus, we point the first feature of solutions to (\ref{q1})
\begin{equation}\label{q3}
u(s)=v(s) \mbox{  \ \ \ \ \ for \ \ } s \in (0,2\pi)\setminus \Xi(u_s).
\end{equation}

From (\ref{q3}) we deduce that  if $s \notin \Xi(v_s)$, then $s\notin
\Xi(u_s)$. Subsequently, we get
$\Xi(v_s) \subset \Xi(u_s)$ which proves the inclusion from (\ref{d5}).
Thus, the isolated elements stay isolated or merge with other elements.
From this we obtain that
$K(v_{s})\geq K(u_{s})$ what ends the proof of line (\ref{d5}).

By properties (\ref{q3}), (\ref{d5}) and the estimate from Lemma
(\ref{provx}), we immediately deduce estimate
(\ref{d3}). In particular, what we gain is a uniform bound in $L_\infty(S)$ on
$\{ \lambda^k_{h,s}\}$. 

The set $\Xi(u_s)$ is defined as the sum of $\bigcup_{l=0}^3 \Xi_l(u_s)$, thus
without loss of generality we can concentrate our attention on one of them,
e.g. on the set $\Xi_2(u_s)$ -- see (\ref{d4}).
From the $\JR$-regularity of $u_s$ set $\Xi_2(u_s)$ is a sum  of closed
intervals, so we take one of them, say,
\begin{equation}\label{q3a}
[a_-,a_+] \subset\Xi_2(u_s)\quad\hbox{ and } u_{s}|_{(a_-,a_+)}=\frac \pi 4.
\end{equation}
Recalling the required regularity of the functions in the $\JR$-class,
%for the fixed interval $[a_-,a_+]$
we find $\epsilon >0$ such that one of the four following possibilities holds:
\begin{equation}\label{q4}
\begin{array}{rcc}
(i) & u_{s}|_{(a_- -\epsilon,a_-)} > \frac \pi 4, &\quad
 u_{s}|_{(a_+, a_++\epsilon)} < \frac \pi 4,
 \\
(ii) & u_{s}|_{(a_- -\epsilon,a_-)} < \frac \pi 4, &\quad
 u_{s}|_{(a_+, a_++\epsilon)} > \frac \pi 4,
 \\
(iii) & u_{s}|_{(a_- -\epsilon,a_-)} > \frac \pi 4, &\quad
 u_{s}|_{(a_+, a_++\epsilon)} > \frac \pi 4,
 \\
(iv) & u_{s}|_{(a_- -\epsilon,a_-)} < \frac \pi 4, &\quad
 u_{s}|_{(a_+, a_++\epsilon)} < \frac \pi 4.
 \end{array}
 \end{equation}
Subsequently, we integrate (\ref{q1}) over $(a_- - \epsilon,a_+ + \epsilon)$
to get
\begin{equation}\label{q5}
\int_{a_- - \epsilon}^{a_+ + \epsilon} u ds - h
(\partial J[u_{s}]|_{a_- - \epsilon}^{a_+ +
\epsilon})=\int_{a_- - \epsilon}^{a_+ + \epsilon} v ds.
\end{equation}

After passing with $\epsilon \to 0^+$, we obtain -- according to the above
four cases (\ref{q1}) -- the following identities
\begin{equation}\label{q6}
\begin{array}{ccc}
(i)&\;\; \int_{a_-}^{a^+} u ds - h \frac \pi 2 = \int_{a_-}^{a^+} v ds&
\mbox{ \ \ (convexity),}
\\
(ii)&\;\; \int_{a_-}^{a^+} u ds + h \frac \pi 2 = \int_{a_-}^{a^+} v ds&
\mbox{ \ \ (concavity),}
\\
(iii) \mbox{ \ and \ } (iv) &\quad \int_{a_-}^{a^+} u ds =
\int_{a_-}^{a^+} v ds&
\mbox{ \ \ (monotonicity).}
\end{array}
\end{equation}

In our present analysis, we essentially use the fact that the energy density
function $J$ is defined by a square. Due to the definition of $J$, see
(\ref{jot}), formula (\ref{q6}) exhausts all the possibilities of the behavior
of $u_s$. For more complex polygons, we would have to discuss more
possible types of facets -- here, there are just  four of them.

We keep considering the interval  $[a_-,a_+]\subset \Xi_2(u_s)$, see
(\ref{q3a}). Let us introduce a set 
\begin{equation}\label{q7}
\Pi=([a_-,a_+]\cap \Xi(v_s))\setminus\left(\Xi_0(u_s) \cup \Xi_1(u_s)\cup
\Xi_3(u_s) \right), 
\end{equation}
then by the properties of sets $\Xi$, we deduce that
\begin{equation}\label{q8}
(u-v)|_{\Pi}= C_h \mbox{ \  is constant}.
\end{equation}
The sign of constant $C_h$ is determined by the geometrical properties of
cases in (\ref{q6}). We have
\begin{equation}\label{q9}
C_h >0 \mbox{ \ \ for \ \ } (i), \qquad
C_h <0 \mbox{ \ \ for \ \ } (ii) \mbox{ \ \ and \ \ }
C_h=0 \mbox{ \ \ for \ \ } (iii) \mbox{ and } (iv).
\end{equation}
Also identity (\ref{q8}) and equation (\ref{q1}) yield
\begin{equation}\label{q10}
 \left.\frac{d}{ds}\partial J[u_{s}]\right|_{\Pi}=\frac{C_h}{h}
\mbox{ \ \ \ \ and \ \ \ \ }
\left.\frac{d}{ds} \partial J[u_{s}]\right|_{(0,2\pi)\setminus \Xi(u_s)}=0.
\end{equation}
Thus, we proved (\ref{d5b}).

Next,  we are going  to study (\ref{d5a}). 
From the analysis of (\ref{q1}), we conclude that
\begin{equation}\label{q11}
||u-v||_{L_1(S)}\leq h \frac{\pi}{2} K(\phi_0).
\end{equation}
Additionally,  from (\ref{d5}) we have also that
$u,v \in W^1_\infty(S)$,
thus simple considerations lead us to the following bound
\begin{equation}\label{q12}
||u-v||_{L_\infty(S)} \leq h^{1/2} C(\phi_0).
\end{equation}

In order to measure the set $\Xi_l(u_s)\setminus \Xi_l(v_s)$ we split it into
two parts 
\begin{equation}\label{q13}
\Xi_l(u_s)\setminus \Xi_l(v_s)=[(\Xi_l(u_s)\setminus \Xi_l(v_s))\cap
  \Xi(v_s)]\;\cup\; 
(\Xi_l(u_s)\setminus \Xi(v_s))=\Pi_1\cup\Pi_2.
\end{equation}

Let us consider $\Pi_1$. On this set we watch the evolution of the intersection
of facets. Thanks to the full information about the  direction of this facet,
we deduce immediately that
\begin{equation}\label{q14}
|\Pi_1|\leq C(\phi_0)h^{1/2}
\end{equation}
The number of possible intersections is controlled by $K(\phi_0)$.

To estimate $\Pi_2$, let us note that this set is a %part %comes directly from 
subset of $\Xi(\lambda_{0,s}+s)$, %for  the initial curve, 
thus in the general case we
can say only 
\begin{equation}\label{q15}
|\Pi_2|\leq V(h),
\end{equation}
where $V(s)\to 0$ as $s \to 0$ and $V$ is determined by the initial datum.
Assuming strict convexity of initial domain we would obtain $V(h)\sim
h^{1/3}$ -- see the example at the end of subsection \ref{sas}.2. \qed

Theorem \ref{thm4xx} is proved.

%---- s 5 -----------
%\section{Convergence of the semi-discretization}\label{scon}
%

\bigskip

Next,  we show that sequences $\{\lambda^k_h\}$ converge to
solutions of the original problem. We will compare solutions given by
Theorem \ref{ist} and Lemma \ref{provx}, in particular, all
assumptions of Theorem \ref{ist} are not required. We follow the
standard procedure which is valid for parabolic operators (see \cite{mra}).

Our next task is to show the following lemma.

\smallskip
%-------lemat x.x rozdz 5.--------
\begin{lemma}\label{lem5xx}
\quad {\sl Let $\Lambda$ and $\{\lambda^k_k\}$ be solutions to
problems (\ref{nasz}) and (\ref{d1}) respectively, then
\begin{equation}\label{c1}
||\Lambda(s,t)-\sum_{k=0}^{[T/h]}\lambda^k_h(s)\chi_{[k,k+1)}(t)||_{
L_1(0,T;L_2(S))} \to 0 \qquad as \quad h\to 0^+.
\end{equation}
If the initial datum fulfills the assumptions of Theorem \ref{ist}, i.e.
$\Lambda_{0,s}\in TV(S)$, then
\begin{equation}\label{c1aa}
||\Lambda(s,t)-\sum_{k=0}^{[T/h]}\lambda^k_h(s)\chi_{[k,k+1)}(t)||_{
L_p(0,T;W^{2-\epsilon}_1(S))} \to 0 \qquad as \quad h\to 0^+
\end{equation}
for any $1<p<\infty$ and $\epsilon >0$.}
\end{lemma}

\smallskip

{\it Proof.} From the properties of solutions to problem
(\ref{nasz}), we know that
$\Lambda_t \in L_2(0,T;L_2(S))$.
It follows that
\begin{equation}\label{c3}
\left|\left|
\frac{\Lambda(s,t)-\Lambda(s,t-h)}{h}-\Lambda_t(s,t)
\right|\right|_{L_1(h,T;L_2(S))} \to 0 \mbox{ \ \ as \  \ }
h \to 0^+.
\end{equation}
For fixed $h>0$ we denote
\begin{equation}\label{c4}
R_h(s,t)=\frac{\Lambda(s,t)-\Lambda(s,t-h)}{h}-\Lambda_t(s,t),
\end{equation}
then the equation (\ref{nasz})${}_1$ can be restated as follows
\begin{equation}\label{c5}
\int_{S}\frac{\Lambda(s,t)-\Lambda(s,t-h)}{h}\pi ds=-\int_{S}
\Om (s,t)
\pi_{s}+R_h(s,t)\pi ds
\end{equation}
for each $\pi$  in $C^\infty(S\times (0,T))$ and each selection $\Om (s,t)$
of multivalued function  $\partial  J[\Lambda_{s}(s,t)+s]$.

We want to compare the above system with the semi-discretization given in
Section \ref{ssem}.
\begin{equation}\label{c6}
\int_{S}\frac{\lambda^k_h(s)-\lambda^{k-1}(s)}{h}\pi ds=
-\int_S \omega(s,t)
%\overline{\partial  J[\lambda^k_{h,s}(s)+s]}
\pi ds
\end{equation}
where $t \in [kh,(k+1)h)$ and $\omega(s,t)$ is any section of $\partial  J[\lambda^k_{h,s}(s)+s]$.

Let us define
\begin{equation}\label{c7}
A^k(s,t)=\Lambda(s,t)-\lambda^k_h(s,t),
\end{equation}
provided $t \in [kh,(k+1)h)$, then from (\ref{c5}) and (\ref{c6}) we deduce
\begin{equation}\label{c8}
\int_S\frac{A^k(s,t)-A^{k-1}(s,t-h)}{h}\pi ds=
-\int_S\left\{ \Om(s,t)
-
\omega(s,t)
\right\} \pi_{s}+R_h(s,t)\pi ds.
\end{equation}

Taking in  (\ref{c8}) as a test function  $A^k(t,s)$,  we get
\begin{equation}\label{c9}
\begin{array}{c}
\int_{(0,2\pi)}|A^k(s,t)|^2 ds=\int_{(0,2\pi)} A^k(s,t)A^{k-1}(s,t-h)ds
\\[8pt]
-h\int_{(0,2\pi)} \left(
\Omega(s,t)-
\omega(s,t)
\right)\left(\Lambda_{s}(s,t)-\lambda^k_{h,s}(s,t)\right)ds
\\[8pt]
+h \int_{(0,2\pi)} R_h(s,t) A^k(s,t)ds,
\end{array}
\end{equation}
but the monotonicity of $\partial J$ implies
\begin{equation}\label{c10}
\int_{(0,2\pi)} \left(
\Omega(s,t)-
\omega(s,t)
\right)
\left(\Lambda_{s}(s,t)-\lambda^k_{s}(s,t)\right)ds \geq 0.
\end{equation}
So, defining
$\alpha^k(t)=||A^k(\cdot,t)||_{L_2(S)},$
by the  Schwarz inequality, we get from (\ref{c4}) the following
inequality
\begin{equation}\label{c12}
\alpha^k(t) \leq \alpha^{k-1}(t-h)+h r^k_h(t) \qquad \mbox{ \ for \ \ }
t \in [kh,(k+1)h),
\end{equation}
where $r^k_h(t)=||R_h(\cdot,t)||_{L_2(0,2\pi)}$.
Thus (\ref{c12}) yields
\begin{equation}\label{c13}
\alpha^k(t) \leq \alpha^0(t-kh)+\sum_{l=1}^k h r^l_h(t-(k-l)h)
\mbox{ \ \  with \ \ }t \in [kh,(k+1)h).
\end{equation}
Integrating (\ref{c13}) over $t \in [kh,(k+1)h)$,
we get
\begin{equation}\label{c14}
\int_{kh}^{(k+1)h} \alpha^k(t)dt\leq \int_0^h \alpha^0(\tau)d \tau +
h \int_0^T ||R_h(\cdot,t)||_{L_2(0,2\pi)} dt
\end{equation}
for $T> (k+1)h$.
Introducing function
$\tilde \alpha(t)=\sum_{l=0}^L \alpha^l(t) \chi_{[lh,(l+1)h)}(t)$
with $L=[T/h]$, from (\ref{c14}) we get
\begin{equation}\label{c16}
\int_0^T \tilde \alpha(t)dt \leq h^{1/2} T
||\Lambda_t||_{L_2(0,T;L_2(0,2\pi))}+ T ||R_h||_{L_1(0,T;L_2(0,2\pi))},
\end{equation}
because the first term of the right-hand-side (RHS) of (\ref{c16})
 is a consequence of the following estimate
\begin{equation}\label{c17}
\begin{array}{l}
\displaystyle
\quad\frac 1h \int_0^h \alpha^0(t)dt \leq \frac 1h \int_0^h
||\int_0^t\Lambda_t(\cdot, \tau)d \tau ||_{L_2(0,2\pi)} dt
\\[8pt]
\displaystyle\leq
\frac 1h (\int_0^h t^2dt)^{1/2} (\int_0^h||\Lambda_t||_{L_2(0,2\pi)}^2
dt)^{1/2} 
\leq C h^{1/2} ||\Lambda_t||_{L_2(0,T;L_2(0,2\pi))} .
\end{array}
\end{equation}

Hence, from (\ref{c16}) and (\ref{c3}) we conclude
$||\tilde \alpha_h||_{L_1(0,T;L_2(0,2\pi))} \to 0 \mbox{ \ \
as \ \ } h \to 0^+$
and we get (\ref{c1}). From the interpolation estimates  and the results
of Theorem \ref{ist}, we deduce that for any $p< \infty$ and
$\epsilon>0$ the convergence (\ref{c1aa}) is valid.
Lemma \ref{lem5xx} is proved.  \qed

%%% --- s 4 ---
\section{Analysis of solutions}\label{sas}

The semi-discretization process proves that the set $\Xi(u^k_{h,s})$ grows with
$k$. One can show that the set $\bigcup_{t\ge0} \Xi(u_{h,s}(\cdot,t))$ may be
estimated from below to show that it  survives the limiting process as $h\to
0$. This may be achieved by the analysis of the semi-discretization procedure,
but this seems tedious. We propose an alternative approach by the construction
of an explicit
solution to (\ref{rn2}) for data in $\vfi_0\in \JR$. By uniqueness
result, see Theorem \ref{ist}, this is the
solution.

We shall assume  in this Section that $\vfi \equiv w_s$ belongs to $\JR$ and
this is the case for the initial data $\vfi_0$ of system (\ref{nasz}).
As a result of the definition of the $\JR$ class we see
\begin{equation}\label{sanr1}
\Xi(\vfi_0) = \bigcup_{l=0}^3 \Xi_l(\vfi_0) =\bigcup_{k=1}^{N_0}
   [\xi_k^-,\xi_k^+],
\end{equation}
where
\begin{equation}\label{sanr8}
\xi_k^-\le \xi_k^+\quad\hbox{and}\quad  \xi_k^+ \le \xi_{k+1}^-, \quad
k=1,\ldots, N_0,
\end{equation}
(with the understanding $\xi^-_{N_0+1} = \xi^-_1+2\pi$). Moreover, each
interval $[\xi_k^-,\xi_k^+]$ is a connected component of one of the sets
$ \Xi_l(\vfi_0)$, $l=0,1,2,3$. We shall also adopt the convention that $0\le
\xi_1^-$ and possibly
$\xi_{N_0}^+ > 2\pi$, but $\xi _{N_0}^+ - 2\pi\le \xi_1^-$.

If $[\xi_k^-,\xi_k^+]$ is one of the connected
components of $ \Xi_l(\vfi)$, then we will  call by a
{\it  facet} the set $F=F_k(\xi_k^-,\xi_k^+)=\{(x,y)\in \bR^2:\ y=w(x), \ x\in
[\xi_k^-,\xi_k^+]\}$. The interval $[\xi_k^-,\xi_k^+]$ will be called the
pre-image of
facet $F_k$. Let us stress  that we admit $\xi_k^-=\xi_k^+$, i.e. a
facet degenerated to a point as well as $\xi^+_{k-1}=\xi^-_k$, i.e. we expect
interaction of facets. We shall see that the generic initial data lead to the
facet
creation (from the degenerate ones) and their interaction. We show
that  facets are formed instantaneously from the data. At this point we
mention that creation of interacting facets leads to additional
difficulties and this process is handled separately.

We will
come up with an explicit formula. Once we check that indeed
this formula yields a solution to equation (\ref{nasz}), we will be assured
that this is  the unique solution we seek.
Subsequently, we shall see that solutions get convexified,
i.e. after some finite time the angle becomes increasing,
hence $w$ becomes convex. Finally, we study
interaction of facets. We will prove that $w(\cdot,t)$ becomes
a minimal solution at the limit time.

It will be also convenient to say that a facet
$F_k(\xi^-_k,\xi^+_k)$, has
{\it zero curvature}, if
$[\xi^-_k,\xi^+_k]$ is a connected component of $\Xi(\vfi)$ and there
exists an open interval $(A,B)$, containing $[\xi^-_k,\xi^+_k]$ such that
$w_s$ is not monotone on any interval $(a,b)$, satisfying
$$
[\xi^-_k,\xi^+_k]\subset (a,b)\subset (A,B).
$$

Furthermore, we say that  a facet $F_k=F_k(\xi^-_k,\xi^+_k)$ is {\it regular}
if
$\xi^-_k< \xi^+_k$. Otherwise, we say that $F_k$ is {\it degenerate}. If
$w_s\in\JR$ is such that the graph of $w$ contains degenerate facets, then we
say that facets are created in solutions to (\ref{nasz}).

Finally, we  say that facets $F_l$, \ldots, $F_{l+r}$ for $r>0$, {\it
  interact} (or are {\it interacting})  if
$F_k\cap F_{k+1}$, $k=l, \ldots, l+r-1$, is a singleton. We call a single
facet $F_k$ {\it non-interacting}, if it is not true that it interacts with any
other facet.

Thus, we have the total of eight combinations, we will treat each case
separately.

\subsection{A  comparison principle}
We are going to establish that solutions to equation (\ref{nasz}) enjoy the expected
comparison principle. This result is interesting for its own sake but
also it is a useful tool analysis. We will apply it to show creation
of interacting facets.

We first recall the basic result (see,
e.g. \cite{smoller}).

\begin{proposition}{\sl Let us suppose that $u_1$, $u_2$ are smooth
    solutions to a strongly parabolic equation
\begin{eqnarray*}
&& u_t =(a(x,u_x))_x \quad\hbox{in }S\times(0,T)
\end{eqnarray*}
and $u_2(x,0) \ge u_1(x,0)$, then  $u_2(x,t) \ge u_1(x,t)$ for all
$t\in(0,T)$. \qed}
\end{proposition}

With this result we may deduce the following comparison principle.

\begin{proposition}{\sl Let us suppose that $\La_1$, $\La_2$ are weak
solutions to (\ref{nasz}) and $\La_1(x,0)\le\La_2(x,0)$, then  $\La_2(x,t)
\ge \La_1(x,t)$ for all $t\in(0,T)$.}
\end{proposition}

\noindent{\it Proof.} Since $\La_1(x,0)\le\La_2(x,0)$, we deduce that
$\La^\ep_1(x,0)\le\La^\ep_2(x,0)$, where $\La^\ep_i$, $i=1,2$ are solutions to
the regularized system (\ref{naszr}). Application of the preceding
result yields
$$
\La^\ep_1(x,t)\le\La^\ep_2(x,t).
$$
Since the point-wise limit exists we conclude that our proposition
holds. \qed

We stress that no information about $\Om_i$, $i=1,2$ is needed in the proof
of the above result.

\subsection{Facet formation}

We shall see below that the evolution of a facet $F_k$
separated from other facets is governed by an ODE for its end-points, see
(\ref{sanr14}) below. In the case
of interacting facets  their evolution is  described by
a system of ODE's (\ref{sanr12}).

As we mentioned we admit facets $F_k$ degenerated to a single point at
the initial instance $t_0=0$. In this case the single ODE (\ref{sanr4}) and
system
ODE (\ref{sanr12}) become singular. While we can resolve satisfactorily the
singularity of the single ODE, the analysis of the system is more
difficult. In fact, we circumvent this problem by using the comparison
principle to show creation of interacting facets.

We shall use the notions and notation introduced above. In addition, in order
to facilitate our construction we shall write
$$
x \mapsto  \al_k(x-s_k)+\tau_k=:l_k(x,s_k,\tau_k) ,
$$
where $\al_k \in {\cal A},$
$s_k\in\Xi(\vfi)$, $\tau_k\in \bR$.

\begin{theorem}\label{sant1}\quad
{\sl Let us assume that  $\vfi_0=w_{0,s}\in\JR$ and $w$ is the unique
 solution to (\ref{nasz}). We also assume that
the  set $\Xi(w_{0,s})= \bigcup_{k=1}^{N_0} [\xi_{k0}^-,\xi_{k0}^+]$ fulfills
 conditions (\ref{sanr1}) and (\ref{sanr8}). Then, there exists a finite
 sequence of time instances $0\le t_0<t_1< \ldots< t_M<\infty$ and a finite
 sequence of continuous functions 
\begin{eqnarray*}
&&\xi_k^\pm:[t_i,t_{i+1}]\to\bR,\qquad i=0,\ldots t_{M-1},\quad
k=1,\ldots N_i,\\
&&\xi_k^\pm:[t_M,\infty)\to\bR,\qquad
k=1,\ldots N_M=4,
\end{eqnarray*}
where $N_0\ge N_1\ge \ldots \ge N_M=4$.

The functions $\xi^-_k(\cdot)$, $\xi^+_k(\cdot)$ satisfying
(\ref{sanr8}) have the following properties:\\
(a) $\xi_k^\pm(0) = \xi_{k0}^\pm$;\\
(b) $0\le \xi_1^-(t)\le\xi_1^+(t)\le  \xi_2^-(t)\le\ldots\le  \xi_N^-(t)\le
\xi_N^+(t)\le \xi_1^- +2\pi$,\ $t\in[t_i,t_{i+1})$;\\
(c) for $t\in[t_i,t_{i+1})$ we have \, $\displaystyle{\Xi(\vfi(\cdot,t))
      = \bigcup_{k=1}^{N_i} [\xi_k^-(t),\xi_k^+(t)]},$
and each interval $[\xi_k^-(t),\xi_k^+(t)]$ is a connected component of one of
the sets $\Xi_l(\vfi(\cdot,t))$, $l=0,1,2,3$.

There exist functions $\tau_k:[t_i,t_{i+1})\to\bR$,\ $i=0,\ldots t_M$,\
$k=1,\ldots N_i,$  and $t_{M+1}=\infty$. They  are such that
the unique  solution to (\ref{nasz}) with initial data  $\vfi(x,0)=\vfi_0(x)$
is given by the following formula for $t\in[t_i,t_{i+1})$, $i=0,\ldots, M$
\begin{equation}\label{sansol}
w(x,t) = \left\{
\begin{array} {ll}
w_0(x) & \hbox{ if }x\in [0,2\pi)\setminus \bigcup_{k=1}^{N_i}
  [\xi_k^-(t),\xi_k^+(t)] \\
l_k(x,\xi_k^+(t_i),\tau_k(t))+w(\xi_k^+(t_i),t_i) & \hbox{ if }x\in
  [\xi_k^-(t),\xi_k^+(t)],\quad
  k=1,\ldots, N_i
\end{array} \right.
\end{equation}
Moreover, $w_x(\cdot,t)$ is well-defined a.e., $\partial w$ defined
by (\ref{depaw}) belongs to $\JR$ and
\newline
\mbox{$\|\partial w(\cdot,t)\|_{\JR}\le
\|\partial w_0\|_{\JR}$}.

In addition, at each time instant $t_i$,
$i=0,\ldots, M$, one of the following happens:

(i) One or more zero-curvature facets disappear, i.e. if one facet disappears
at $t_i$, then
$$
\xi_{k_0-1}^+(t)\le \xi_{k_0}^- (t) <  \xi_{k_0}^+ (t) \le
\xi_{k_0+1}^-(t),\quad \hbox{for } t_i<t<t_{i+1}
$$
and
$$
\lim_{t\to t_{i+1}^-} \xi_{k_0-1}^- (t) =  \xi_{l_0}^- (t_{i+1}),\qquad
\lim_{t\to t_{i+1}^-} \xi_{k_0+1}^- (t) =  \xi_{l_0}^+ (t_{i+1}),
$$
where $[\xi_{l_0}^-(t_{i+1}),  \xi_{l_0}^+ (t_{i+1})]$ is a subset of a
connected component of $\Xi_l(\vfi(t_{i+1}))$, as a result $N_{i+1} < N_i$.

(ii) One pair or more pairs of facets begin to interact, i.e.
$\xi_{k-1}^+(t) < \xi_{k}^-(t)$ for $t_i<t<t_{i+1}$ and
$$
\lim_{t\to t_{i+1}^-} \xi_{k-1}^+ (t) = \xi_{k-1}^+ (t_{i+1}) = \xi_{k}^-
(t_{i+1}) = \lim_{t\to t_{i+1}^-} \xi_{k}^- (t).
$$
}

\end{theorem}
\noindent

The proof is achieved in a number of steps. Its major parts are
separated as Lemmas. We start with constructing the $\xi_k^\pm$'s. We
first consider non-interaction during creation of facets, i.e.
\begin{equation}\label{sanr13}
\hbox{if }\xi^-_k=\xi^+_k,\quad\hbox{then}\quad \xi^+_{k-1}<\xi^-_k
\hbox{ and }\xi^+_{k}<\xi^-_{k+1}.
\end{equation}
However, the lemma below is valid without this restriction.

\begin{lemma}\label{sanl1}\quad{\sl
Let us suppose that $w_s =\vfi\in\JR$  and $[\xi_k^-,\xi_k^+]$ is
a connected component of $\Xi_l(\vfi)$ and $s_k$ is its member. We assume
that $F_k(\xi_k^-,\xi_k^+)$ is not a zero curvature facet. \\
(a) If $\xi_k^+ < \xi_{k+1}^-$, then for sufficiently small $\tau_k$ of a
proper sign, there exist $\xi_k^\pm(\tau_k)$ such that
\begin{equation}\label{sanr2}
w(\xi_k^\pm(\tau_k)) = l_k(\xi_k^\pm(\tau_k), \xi_k^\pm, \tau_k)+
w(\xi_k^\pm)\qquad\hbox{and }\xi_k^\pm(0) = \xi_k^\pm.
\end{equation}
Moreover, the functions $\tau_k\mapsto \xi_k^\pm(\tau_k)$ are Lipschitz
continuous, provided that $w_s(\xi^\pm_k) \ne \al_k$. Otherwise,
$\xi_k^\pm(\tau_k)$ are locally Lipschitz continuous. In  addition,
\begin{equation}\label{sanr3}
\frac{d \xi_k^+}{d\tau_k}(\tau_k) = \frac 1{w_s( \xi_k^+) -\al_k},
\quad\frac{d \xi_k^-}{d\tau_k}(\tau_k) = \frac 1{w_s( \xi_k^-) -\al_k}
\qquad \hbox{for } a.e. \ |\tau_k| \in [0,\ep).
\end{equation}
(b)  If $\xi_l^+ < \xi_{l+1}^-\le  \xi_{l+1}^+ =
\xi_{l+2}^-\le\xi_{l+2}^+= \xi_{l+3}^-
\ldots \le
\xi_{l+r}^+  < \xi_{l+r+1}^- $, (in particular we admit
$\xi_{1}^-=\xi_{N_0}^+ -2\pi$), then
\begin{equation}\label{sanr9}
w(\xi_{k-1}^+) +
l_{k-1}(\xi_{k-1}^+(\tau_{k-1},\tau_{k})),\xi_{k-1}^+,\tau_{k-1}) =
w(\xi_{k}^+) + l_{k}(\xi_{k}^+(\tau_{k},\tau_{k+1})),\xi_{k}^+,\tau_{k})
\end{equation}
for $k=l+1,\ldots,l+r$.
Moreover, the functions $(\tau_k,\tau_{k+1})\mapsto
\xi_k^\pm(\tau_k,\tau_{k+1})$, $k=l+1,\ldots,l+r-1$ are Lipschitz continuous.
}
\end{lemma}

\bigskip\noindent{\it Proof.}
Before proceeding to the formal proof we will explain the situation by
drawing a picture (where the subscript $k$ is suppressed). The graph of
$w(\cdot)$ and the line containing $F(\xi^-,\xi^+)$ moved vertically by $\tau$
intersect at $x=\xi^-(\tau)$ and at  $x=\xi^+(\tau)$.

\bigskip
%\centerline{
%{\epsfxsize=8cm \epsfbox{mr-fig.eps}}}

\bigskip

(a) Since $F_k$ is not of zero curvature then by
the fact that $\vfi\in\JR$ it follows that $w$ in a neighborhood of
$[\xi_k^-,\xi_k^+]$ is either convex or concave. Let us consider the
case of $w$ being convex on $(a,b)\supset [\xi_k^-,\xi_k^+]$, the
other case is similar. By convexity, any chord is above the graph of
$w$. Thus,  the line $l_k(\cdot,\xi^+_k,\tau_k)+w(\xi_k^+)$ for
sufficiently small $\tau_k>0$ intersects the graph of $w$ at exactly
two points, i.e. for $\tau_k>0$ equation  (\ref{sanr2}) has exactly
two solutions. One of them, which is greater than $\xi^+_k$ is called
$\xi^+_k(\tau_k)$, the other one, smaller than $\xi_k^-$ is dubbed
$\xi_k^-(\tau_k)$. The
function
\begin{equation}\label{sanr10}
x\mapsto w(x) - l_k(x,\xi^+_k,\tau_k)-w(\xi_k^+)=: F^+_k(x)
\end{equation}
is increasing for $x\in[ \xi_k^+,b)$ and this interval is maximal with
this property, while the function
\begin{equation}\label{sanr101}
x\mapsto w(x) - l_k(x,\xi^-_k,\tau_k)-w(\xi_k^-)=: F^-_k(x)
\end{equation}
and decreasing  for $x\in(a, \xi_k^-]$ and again this interval is maximal with
this property. One can see this by taking the
derivative of  (\ref{sanr10}) and (\ref{sanr101}), because we have
$$
\frac d{dx}(w(x) - l_k(x,\xi^+_k,\tau_k)) = w'(x)-\al_k\ge
w'(\xi^+_k)-\al_k > 0\qquad\hbox{for } a.e.\ x\in[\xi^+_k,b)
$$
and
$$
\frac d{dx}(w(x) - l_k(x,\xi^-_k,\tau_k)) = w'(x)-\al_k\le
w'(\xi^-_k)-\al_k < 0\qquad\hbox{for } a.e.\ x\in(a,\xi^-_k].
$$
Thus, the function $[\xi_k^+, b) \in x \mapsto F_k^+ (x)$
  (resp. $(a,\xi_k^-] \mapsto F_k^- (x)$  has a continuous inverse. As
a result, for any $\tau_k$ belonging to $[0,\delta) \subset F_k^+
  ([\xi_k^+, b)) \cap F_k^-((a,\xi_k^-])$, $\delta>0$,
we may set $\xi_k^+(\tau_k) = (F^+_k)^{-1}(\tau_k) $ and
$\xi_k^-(\tau_k) = (F^-_k)^{-1}(\tau_k) $. Moreover,
$$
\frac {d\xi_k^\pm}{d\tau_k}(\tau_k)= \frac 1{w_s(\xi^\pm_k(\tau_k))-\al_k},
\quad a.e.
$$
This formula combined with monotonicity of $w_s$ yields,
\begin{equation}\label{sanr11}
\frac 1{\al_k-w^-_s(b)}\le \frac {d\xi_k^+}{d\tau_k}(\tau_k)
\le \frac 1{w^+_s(\xi_k^+(\tau_k)) -\al_k},
\quad
\frac 1{w^+_s(a)-\al_k}\le\left| \frac {d\xi_k^-}{d\tau_k}(\tau_k)\right|
\le \frac 1{w^-_s(\xi_k^-(\tau_k)) -\al_k}
\end {equation}
for a.e $\tau_k$. If $w_s(\xi^+_k)\ne \al_k$, then it follows that
$\xi^+_k(\cdot)$ is Lipschitz continuous on
$[0,\ep]%xi^+_k,\xi^+_k+\ep]
$,
for some $\ep>0$. A similar statement is valid for $\xi^-_k(\cdot)$.

(b) Functions $\xi^+_{l}(\tau_{l},\tau_{l+1})$, \ldots,
$\xi^+_{l+r-1}(\tau_{l+r-1},\tau_{l+r})$ are defined as unique solutions to
the decoupled system of linear equations (\ref{sanr9}) for any given
$\tau_{l},$ $\ldots$, $\tau_{l+r}$.
This is indeed possible because $\al_k\ne\al_{k+1}$. The solution
$\xi_k^+$ depends linearly upon $\tau_k$, $\tau_{k+1}$. Subsequently, we set
$\xi^-_{k+1}:=\xi^+_k(\tau_{k},\tau_{k+1})$, $k=l ,\ldots,l+r-1$.
\qed

\bigskip\noindent{\bf Remark.} In the case (a) the derivatives $\frac
d{d\tau_k}\xi_k^\pm$ are never zero. They may converge to infinity at $t=t_i$,
as well as at $t=t_i+t^*$,
if at that time instance $w_s(\xi_k^\pm) = \al_k$.

\bigskip
The lemma above expressed the evolution of the pre-images of facets in terms of
$\tau_k$, i.e. the amount of vertical shift of the line
$l_k(\cdot,\xi^-,w(\xi^-))$. However, in order to render (\ref{sansol})
meaningful, we have to figure out the time dependence of $\tau_k$. At the same
time we have to construct $\Om$. We begin with an explicit case.

\begin{lemma}\label{sanl2}
{\sl Let us suppose that $F_k(\xi_k^-,\xi^+_k)$ is neither of zero-curvature
    nor interacting and it may be degenerate. Then there exist $\Om_k^-$,
    $\Om_k^+\in\partial J(\al_k)$ and a unique solution
    $\tau_k:[t_*,t_*+T_{max})\to\bR$ to the equation
\begin{equation}\label{sanr4}
\frac{d\tau_k}{dt} = \frac{\Om_k^+ -\Om_k^-}
{\xi_k^+(\tau_k)- \xi_k^-(\tau_k)},\qquad \tau_k(t_*) =0.
\end{equation}
They are such that the function
\begin{equation}\label{sanr5}
\Om(x,t) = \frac{\Om^+_k-\Om^-_k}{\xi^+_k(t)-\xi^-_k(t)}(x-\xi^-_k(t))
+ \Om^-_k
\end{equation}
and $w$ defined  by (\ref{sansol}) satisfy}
$$
\frac{\partial w}{\partial t}(s,t) = \frac{\partial \Om}{\partial s}(s,t)
\qquad\hbox{\sl for }t\in [t_*,t_*+T_{max}),\quad s\in
  (\xi_k^-(\tau_k(t)),\xi_k^+(\tau_k(t)).
$$
\end{lemma}
\bigskip\noindent{\it Proof.} %We first handle the case
%facets  $F_k(\xi_k^-,\xi_k^+)$
%which do not have zero curvature. We have to consider two sub-cases:
%(a) $\xi_k^+< \xi_{k+1}^-$; (b)  $\xi_k^+ = \xi_{k+1}^-$. In the first case
%we set
The non-interaction assumption implies that
$$
\al_{k-1} \equiv\al_k-\De\al < w_{0,s}^-(\xi^-)\le \al_k\le
w_{0,s}^+(\xi^+)<\al_{k+1} \equiv\al_k+\De\al
$$
or
$$
\al_{k-1} > w_{0,s}^-(\xi^-)\ge \al_k\ge
w_{0,s}^+(\xi^+)>\al_{k+1}.
$$
Keeping this in mind we set
\begin{equation}\label{rnpop0}
\Om_k^+ =\lim_{x\to (\xi_k^+(t))^+} \frac{\partial J}{\partial
  \vfi}(w_{0,x}(x))\qquad
 \Om_k^- =\lim_{x\to (\xi_k^-(t))^-} \frac{\partial J}{\partial
  \vfi}(w_{0,x}(x)).
\end{equation}
Of course $ \Om_k^-,$ $\Om_k^+\in \partial J(\al_k)$. We notice that both
quantities are well-defined for regular as well as degenerate facets.

%If (b) holds, then we set
%\begin{equation}
%\{\Om(t,\xi_k^+(t))\} = \partial J(\al_k)\cap  \partial J(\al_{k+1}).
%\end{equation}
%We proceed in a similar manner to define $\Om(t,\xi_k^-(t))$ as an
%element of $\partial J(\al_k)\cap  \partial J(\al_{k-1})$. This
%intersections contain exactly one point for $t_i>0$.

%In order to construct solutions to (\ref{nasz}) we have to find the time
%dependence of $\tau_k$. For this purpose we consider the following ODE
%for $\tau_k$,
%
%This definition will become more clear later in the proof of Theorem
%\ref{sant1}, see \ref{bullet} below.
%[KIEDY LETTER??]
Now, we turn our attention to equation (\ref{sanr4}), we notice that this
equation states that the time derivative
of $\tau_k$ equals the slope of the straight line passing trough the points
$(\xi_k^-,\Om(\xi_k^-))$ and $(\xi_k^+,\Om(\xi_k^+))$ . This line provides a
section of $\partial J$, necessary to construct solutions to (\ref{nasz}).

The numerator of (\ref{sanr4}) is constant and if
$\xi^+_k(\cdot)$, $\xi^-_k(\cdot)$ are Lipschitz continuous and
$\xi^+_k(\tau_k) > \xi^-_k(\tau_k)$ for all the values of $\tau_k$,
then (\ref{sanr4}) has a unique solution. If however, $\xi^+_k(0)=\xi^-_k(0)$,
then (\ref{sanr4}) is singular and this equation requires special attention. A
similar situation arises when $w_s(\xi^\pm_k) = \al_k$.
Fortunately, due to a simple structure of (\ref{sanr4}) we may resolve these
issues.

The ODE (\ref{sanr4}) governing the behavior of  a non-interacting facet $F_k$
is obtained by taking the time derivative of (\ref{sanr2}),
\begin{equation}\label{sanr14}
\frac d{dt}(w(\xi_k^+(t)) -\al_k\xi_k^+(t)) =
\frac d{dt}\tau_k(t),\qquad
\frac d{dt}(w(\xi_k^-(t)) -\al_k\xi_k^-(t)) =
\frac d{dt}\tau_k(t)
.
\end{equation}
In reality, we do not assume that $w$ is differentiable everywhere, but its
one-sided derivatives do exist at each point. Due to monotonicity  of
$\xi^\pm_k$ the one-sided derivatives suffice in the formula above.

%\begin{lemma}\label{sanl2}\quad
%{\sl Let us suppose that facet $F_k(\xi^-_k,\xi^+_k)$ has non-zero
%  curvature and it is non-interacting at $t_0$ and
%  $\xi^-_k\le\xi^+_k$, and the equality is admitted. Then, there
%  exists $\tau_k$ a unique solution to (\ref{sanr4}).}
% two continuous functions $\xi^\pm_k:[t_0,t_1)\to \bR$ such that they
% are locally Lipschitz continuous on $(t_0,t_1)$ and satisfy
% (\ref{sanr14}) a.e.
%\end{lemma}
%\bigskip\noindent{\it Proof.}
By the definition of $\xi^\pm_k$ we rewrite (\ref{sanr4}) as follows
$$
((F_k^+ )^{-1}(\tau_k) - (F_k^- )^{-1}(\tau_k)) \frac{d\tau_k}{dt} =
\De\Om_k.
$$
Here, due to  the definition of $J$ and (\ref{rnpop0}), we  have
$$
\De\Om =\De\Om_k =  \Om^+ - \Om^-= \frac\pi 2.
$$
Since $\xi^+_k(\tau_k) > \xi^-_k(\tau_k)$ as long as $\tau_k\neq 0$,
then we deduce that $G$,
the primitive function of $(F_k^+ )^{-1}(\tau_k)-(F_k^-)^{-1}(\tau_k)$
such that $G(0)=0$, is strictly increasing. Thus (\ref{sanr4})
takes the form
$$
\frac d{dt}\left( G(\tau_k) \right) = \De\Om
$$
or $G(\tau_k)= \De\Om t$. As a result function $\tau_k$ is given uniquely
by the formula
$$
\tau_k(t) = G^{-1}(\De\Om t)
$$
and $\tau_k(0) =0$.

If we now set $\Om$ by formula (\ref{sanr5}), then by the convexity of
the set $\partial
J(\al_k)$, we conclude  that $\Om(x,t)\in \partial J(\al_k)$. Moreover, for $w$
defined by (\ref{sansol}), the  following equality holds by the
definition of $\Om$ and $\tau_k$,
$$
\frac{\partial w}{\partial t}(x,t) = \frac{d\tau_k}{dt}
= \frac{\Om^+_k - \Om^-_k}{\xi^+_k(t)- \xi^-_k(t)} =
\frac{\partial\Om} {\partial x}(x,t)
$$
for $t\in[t_*,t_*+T_{max})$, $x\in (\xi^-_k(t), \xi^+_k(t))$.

We note that $\Om$, which we so far constructed, belongs to
$W^{1}_{1}([0,2\pi))$ for each $t>t_*$, if however the facet does not
degenerate, then $\Om(\cdot,t_*) \in W^{1}_{1}([0,2\pi))$ too.

\qed
%If we take into account the definition of $F^\pm_k$, then we can set
%$$
%\eta^+ = F^+_k(\xi^+_k)\quad\hbox{and}\quad
%\eta^- = F^-_k(\xi^-_k).
%$$
%These two new variables are well-defined.
%
%We want to write (\ref{sanr14}) using $\eta^\pm$ and the formula for
%$\frac{d\tau_k}{dt}$. However, in order to proceed we have to make
%sense out of the numerator of the fraction in (\ref{sanr4}) at
%$t_0$. We know that is constant equal to $\De\Om_k$ for $t>t_0$, thus
%we define by continuity the value of this numerator at $t_0$ to be
%$\De\Om_k$. Hence (\ref{sanr14}) takes the form of a system for
%$\eta^+$, $\eta^-$,
%$$
%\frac d{dt}\eta^+ =
%\frac{\De\Om}{(F_k^+)^{-1}(\eta^+)-(F_k^-)^{-1}(\eta^-)},\qquad
%\frac d{dt}\eta^- =
%\frac{\De\Om}{(F_k^+)^{-1}(\eta^+)-(F_k^-)^{-1}(\eta^-)}.
%$$
%We immediately conclude from this system that $\eta^+-\eta^-=const$
%and furthermore $\eta^+(t) =\eta^-(t)$.  Hence, we can eliminate
%$\eta^-$ and consider only
%\begin{equation}\label{sanr15}
%\frac d{dt}\eta^+ =
%\frac{\De\Om}{(F_k^+)^{-1}(\eta^+)-(F_k^-)^{-1}(-\eta^+)},\qquad
%\eta^+(0)=F_k^+(\xi^+_k).
%\end{equation}
%The denominator is positive and continuous, hence it the derivative of
%an increasing function, denoted by $G(\eta^+)$. Thus, (\ref{sanr15})
%becomes
%$$
%G(\eta^+(t)) - G(\eta^+(t_0)) = \De\Om(t-t_0).
%$$
%Our claim follows.

We can infer the following observation from Lemma \ref{sanl1} and
(\ref{sanr4}).

\begin{corollary}\label{sanc1}\quad{\sl
Let us suppose that $w_s$ is increasing (resp. decreasing) in a neighborhood
of the pre-image $[\xi^-_k,\xi_k^+]$ of a non-interacting facet. Then,
there exists a positive $\de$, such that for $t\in[t_k,t_k+T)$:\\
(a) if $\xi_k^+ < \xi_{k+1}^-$, then $\frac d{dt}\xi_k^+(\tau_k(t))\ge \de>0$
a.e. (resp. $\frac d{dt}\xi_k^+(\tau_k(t))\le \de<0$ a.e.).\\
(b) if $\xi_{k-1}^+ < \xi_{k}^-$, then $\frac d{dt}\xi_k^-(\tau_k(t))\le \de<0$
a.e. (resp. $\frac d{dt}\xi_k^-(\tau_k(t))\ge\de>0$ a.e.).
}
\end{corollary}

\medskip\noindent{\it Proof.} The chain formula yields
$\frac d{dt}\xi_k^+ = \frac {d\xi_k^+}{d\tau_k} \frac {d\tau_k}
      {dt}\quad a.e.$
In the case (a), by the geometry of the problem, we deduce that
$\frac {d\xi_k^+}{d\tau_k}>0$ (see  (\ref{sanr3})) as well as  $\frac
{d\tau_k} {dt}>0$ (see (\ref{sanr4})). Moreover, formulas  (\ref{sanr3}) and
(\ref{sanr4}) imply that none of the factors may vanish, in fact they are
separated from zero.

The remaining cases are handled in the same way. \qed

\medskip
We shall state a result corresponding to Lemma \ref{sanl2} for a set of
interacting facets. It will be somewhat more tedious.

\begin{lemma}\label{sanl3}\quad
{\sl Let us suppose that non-degenerate facets $F_l$, \ldots, $F_{l+r}$, $r>0$
interact,
while $\xi^+_{l-1} < \xi^-_l$ and $\xi^+_{l+r}<\xi^-_{l+r+1}$. Then, there
exist continuous functions $\xi^\pm_k:[t_*,t_*+T)\to \bR$, $k=l,$ \ldots,
${l+r}$, such that they are locally Lipschitz continuous on
$(t_*,t_*+T)$ satisfying (\ref{sanr12}) below and there are $C^1$ functions
  $\tau_k:[t_*,t_*+T)\to \bR$, $k=l,$ \ldots, ${l+r}$, and
$\Om(\cdot,t)\in W^1_1(\xi^-_l(t),\xi^+_{l+r}(t))$. They are
all such that $w$ defined  by (\ref{sansol}) satisfies}
\begin{equation}\label{rnpop3}
\frac{\partial w}{\partial t}(s,t) = \frac{\partial \Om}{\partial s}(s,t)
\qquad\hbox{\sl for }t\in [t_*,t_*+T_{max}),\quad s\in
  (\xi_l^-(t),\xi_{l+r}^+(t)).
\end{equation}
\end{lemma}

\smallskip\noindent{\bf Remark.} The above Lemma includes the case when the
set $S\setminus\Xi(w_{0,s}) $ consists of a single component.

\bigskip\noindent{\it Proof.} By our assumption the pairs of facets
$F_{l-1}$, $F_l$ and $F_{l+r}$, $F_{l+r+1}$ do not interact. Thus, the
evolution of the end points $\xi^-_l$ and $\xi^+_{l+r}$ is determined
as for a single non-interacting facet. This remain applicable, unless
$\Xi(w_s(\cdot,t_i))=[0,2\pi)$. We proceed as in Lemma \ref{sanl2}, but we have
to  determine   $\xi_{l+i}^\pm$, $\tau_{l+i}$, $i=1,$ \ldots, $r$  and $\Om$
simultaneously. We keep in mind that $\xi^-_{l+i}=\xi^+_{l+i-1}$, $i=1,$
\ldots,  $r$.
In order to obtain their time evolution,  we differentiate (\ref{sanr9})
with respect to  time. This yields,
\begin{equation}\label{rnpop2}
\al_k\dot\xi^+_k +\dot\tau_k = \al_{k+1}\dot\xi^+_{k+1} +\dot\tau_{k+1}.
\end{equation}
The equation for $\tau_k$ should be similar to (\ref{sanr4}), if so we have to
select $\Om_{l+i}^\pm$, $i=1,$ \ldots, $r$. We define $\Om_l^-$ and
$\Om_{l+r}^+$ as in  (\ref{rnpop0}), i.e.
\begin{equation}\label{rnpop1}
\Om_{l+r}^+ =\lim_{x\to (\xi_{l+r}^+(t))^+} \frac{\partial J}{\partial
  \vfi}(w_{0,x}(x)),\qquad
 \Om_l^- =\lim_{x\to (\xi_k^-(t))^-} \frac{\partial J}{\partial
  \vfi}(w_{0,x}(x)).
\end{equation}
We have to define the remaining $\Om^\pm_k$'s while keeping in mind
$\Om_k^+ = \Om_{k+1}^-$. By the properties of derivative
$\frac{\partial J}{\partial \vfi}(\vfi) $ and the subdifferential
$\partial J(\al_k) $ the number $\Om^-_l$ is one endpoint of the interval
$\partial J(\al_k) $, thus we inductively define $\Om_k^+$ as follows,
\begin{equation*}
\Om_{k+1}^+ =\left\{
\begin{array}{ll}
\Om_k^+, & \hbox{if the facet } F_k \hbox{ has zero curvature}\\
\hbox{the other endpoint of the interval } \partial J(\al_k), &
\hbox{otherwise.}
\end{array}
\right.
\end{equation*}
We have to check that $\Om_{l+r}^+$ defined in this way agrees with
(\ref{rnpop1})${}_2$. We prove this by induction with respect to $r$,
the number of interacting facets. If $r=1$, then the claim follows
from the preceding 
considerations. Let us suppose validity of the claim for some $r\ge1$, we will
show it for $r+1$. Let us suppose that $w_0$ corresponds to a group of $r+1$
interacting facets satisfying the assumptions of the Lemma. We consider such a
mollification $w^\ep_0 $ of   $w_0 $  in a neighborhood of $\xi^+_{l+r} =
\xi^-_{l+r+1}$ that $w^\ep_0 = w_0$ for $x$ satisfying $|x-
\xi^+_{l+r}|\ge\ep$ and $w^\ep_0 $ is smooth. Moreover, we require that
$w_{0,s}$ and $w^\ep_{0,s} $ are simultaneously increasing or decreasing. Thus
the facets corresponding to $w^\ep_0$ are $F_l$, \ldots, $\tilde F_{l+r}$,
$\tilde F_{l+r+1}$. We notice that facet  $\tilde F_{l+r+i}$ is of zero
curvature iff facet  $F_{l+r+i}$ is of zero curvature, $i=0,1$. Moreover,
facets $\tilde F_{l+r}$,
$\tilde F_{l+r+1}$ do not interact. By the inductive assumption
$\tilde\Om^+_{r+l}=\Om^+_{r+l}$  is equal to $\Om^-_{r+l+1}$. At the same time
$\tilde\Om^+_{r+l} = \tilde\Om^-_{r+l+1}$ is determined from
$\tilde\Om^+_{r+l+1}= \Om^+_{r+l+1}$ and $w_{0,s}$ as in Lemma
\ref{sanl2}. The two ways of course coincide, due to formulae
(\ref{rnpop0}). Our claim follows.

We now write equations for $\tau_k$, $k=l$, \ldots, $l+r$, they are as
(\ref{sanr4}),
\begin{equation}\label{uta}
\frac{d\tau_k}{dt} = \frac{\Om_k^+ -\Om_k^-}
{\xi_k^+(\tau_k)- \xi_k^-(\tau_k)},\qquad \tau_k(t_*) =0\qquad
\mbox{ \ \ for \ \ } k=l, \ldots,
l+r.
\end{equation}
Since we do not admit degenerate facets, these equations are not singular.
We combine them with (\ref{rnpop0}) and after writing
$\eta=(\xi_{l}^+,\ldots,\xi^+_{l+r-1})$, we arrive at
\begin{equation}\label{sanr12}
A\dot\eta = B(\eta),
\end{equation}
where
$$
A =\left[
\begin{matrix}
\al_{l} & -\al_{l+1} & \ddots & 0  \\
0 & \al_{l+1} & -\al_{l+2} & 0\\
\ldots &    \ddots & \vdots & -\al_{l+r-1} \\
 \ddots&          0  &  0&\al_{l+r-1}
\end{matrix}\right],
$$
$$
B(\eta)_k = -\frac{\Om_k^+ - \Om_k^-}{\eta_k^+ - \eta_{k-1}^+}
+\frac{\Om_{k+1}^+ - \Om_{k+1}^-}{\eta_{k+1}^+ - \eta_{k}^+}, \qquad
k=l+1,\ldots,l+r-2,
$$
$$
B(\eta)_{l+r-1} = \al_{l+r}\frac d{dt}\xi^+_{l+r}
-\frac{\Om_{l+r-1}^+ - \Om_{l+r-1}^-}{\eta_{l+r-1}^+ - \eta_{l+r-2}^+}
+\frac{\Om_{l+r}^+ - \Om_{l+r}^-}{\xi_{l+r}^+ - \eta_{l+r-1}^+}
.
$$
Under our assumptions, there is a separate equation for $\frac
d{dt}\xi^+_{l+r}$ i.e. (\ref{sanr4}). Due to the assumption of absence
of degenerate interacting
facets, this system is uniquely solvable on $[t_*, t_*+T)$.

We have to define $\Om$, it will be a continuous piece-wise linear function,
\begin{equation}\label{omega}
\Om(x,t) =
\frac{\Om_{l+i}^+ - \Om_{l+i}^-}{\xi_{l+i}^+(t) - {\xi_{l+i}^-(t)}}
(x-\xi_{l+i}^-(t)) + \Om_{l+i}^-.
\end{equation}
Moreover, $w$ and $\Om$ satisfy  (\ref{rnpop3}). \qed

We claim in theorem \ref{sant1}that the number of facets decreases in time. The
result below explains that certain phenomena are forbidden. Namely, no facet
with non-zero curvature  may degenerate.

\begin{proposition}\label{sanp1}{\sl
In any group of interacting facets $F_k$, $k=l,\ldots,l+r$, $r>0$ only a facet
with zero curvature may degenerate.
%It is not possible that an interacting facet $F_k$ with non-zero curvature
%becomes degenerate after some time of evolution.
}
\end{proposition}\bigskip\noindent{\it Proof.} Let us suppose  that
$F_l$, $\ldots$, $F_{l+r}$, $r>0$ is a maximal group of interacting facets with
non-zero curvature. For the sake of definiteness, we will proceed while
assuming that $w_s$ is increasing on $(a,b)\supset [\xi_l^-,\xi^+_{l+r}]$.

%%We sIn order to verify it we look at (\ref{rnuz0}) and
%%notive that for $t$ close to $T$

\noindent{\it Step 1.} Let us observe that for a facet $F_k$ to
disappear, it is necessary, (but not sufficient) that one of neighboring facets
moves upward faster than $F_k$, i.e. either
$V_{k+1}=\frac{d\tau_{k+1}}{dt}>\frac{d\tau_{k}}{dt} =V_k $ or
 $V_{k-1}=\frac{d\tau_{k-1}}{dt}>\frac{d\tau_{k}}{dt} =V_k $. Indeed,
the position of $F_k$ is defined by the intersection of the lines
containing $F_k$, $F_{k+1}$ moved vertically by $\tau_k$ and
respectively by  $\tau_{k+1}$ and the intersection of
lines containing $F_k$, $F_{k-1}$ moved vertically by $\tau_{k}$ and
respectively by  $\tau_{k-1}$. Thus,
if the lines  containing $F_{k+1}$ and  $F_{k-1}$ are moved up so
much that their intersection is above the line  containing $F_k$ moved
vertically by $\tau_{k}$, then facet $F_k$ is going to disappear. This
situation may occur only if $V_{k+1}>V_k$ or $V_{k-1}>V_k$.

\noindent{\it Step 2.} Let us suppose that facets $F_k$, $F_{k-1}$
interact, hence by (\ref{sanr12})
\begin{equation}\label{rnuz0}
\al_k\dot\xi_k - \al_{k-1}\dot\xi_{k-1} = \dot\tau_{k-1} -\dot\tau_k.
\end{equation}
By the monotonicity assumption on $w_s$ we notice that
$\dot\tau_{k-1}$ and $\dot\tau_k$ are positive.
If the length of $F_k$, which is equal to $\xi_k-\xi_{k-1}$, stays
bounded on $[t_*,t_*+T)$ while the length
of $F_{k-1}$ vanishes at $t=t_*+T$, then in a neighborhood of $t_*+T$ we have
$\dot\tau_{k-1} -\dot\tau_k<0$.
Thus, by (\ref{rnuz0}) we can see that
$$
\al_{k-1} (\dot\xi_k-\dot\xi_{k-1}) + (\al_k- \al_{k-1}) \dot\xi_k  <0
$$
and by (\ref{uta}) the left-hand-side (LHS) converges to $-\infty$ when $t$
tends to 
$t_*+T$. Since $\dot\xi_k-\dot\xi_{k-1}$ must be bounded from above, we
deduce that $\dot\xi_k<0$ for $t$ close to $t_*+T$.

\noindent{\it Step 3.} Since always $\dot \xi_{l-1}^-<0$ and $\dot
\xi_{l+r+1}^+>0$ (unless $\Xi(\vfi) = [0,2\pi)$), we conclude that not
all of the facets vanish simultaneously at $t=t_*+T$. As a result we may
assume the length $\ell(F_{l-1})$ of $F_{l-1}$ is greater than $d>0$
on $[t_i,T)$. Thus, we conclude by step 1, that for $t$ close to $t_*+T$
we have $V_l>V_{l-1}$. By induction we obtain that
\begin{equation}\label{rnuz1}
V_{k+1} >V_k, \qquad k=l,\ldots, j+r-1.
\end {equation}

We notice that we have the following possibilities for facet
$F_{l+r}$: (a) there is an adjacent zero-curvature facet $F_{l+r+1}$;
(b) $\xi^+_{l+r+1}$ is defined as $(F_{l+r+1}^+)^{-1}(\tau_{l+r+1})$
(see the proof of Lemma \ref{sanl1}). In case (a) we can see that
$\tau_{l+r+1}=0$ while in (b) $\tau_{l+r+1}>0$.

The condition  (\ref{rnuz1}) combined with (\ref{uta}) implies that
$$
\xi^+_l -\xi^-_l > \ldots > \xi^+_{r+l+1} -\xi^-_{r+l+1}.
$$
Hence, the endpoints of $F_k$, $  k=l,\ldots, j+r-1$ converge to a
common limit $p$. But by step 2
$$
 \xi^+_{r+l+1}(t) >  \xi^+_{r+l+1}(t_*) >  \xi^+_{l+1}(t_*) >
 \xi^+_{l+1}(t).
$$
This is a contradiction, our claim follows. \qed

This observation shows that the initial time $t_0=0$ is special. If
the data are poor 
from the view-point of dynamics, but still acceptable, then they get
immediately regularized. That is all non-zero curvature degenerate facet
become regular.

We are now ready for the proof of the main result.

\bigskip\noindent{\it Proof of Theorem \ref{sant1}.} {\sc Part A.}
We start with data free from degenerate interacting facets.
We set $t_0=0$, we    have to define time instance $t_i$, $i=1,\ldots,M$
postulated by the theorem. We shall proceed iteratively.

It follows from Proposition \ref{sanp1}, that degenerate, non-zero curvature
facets are possible only at $t=0$, i.e. at the initial time instance.
%, but first assume that
%$t_i>0$ is given, then we construct $t_{i+1}$, because $t_0$ requires
%additional considerations.

Let us suppose %then,
that %$\xi_k^-(t_i)< \xi_k^+(t_i)$ and
$[\xi_k^-(t_i), \xi_k^+(t_i)]$ is a connected component of $\Xi_l(w_s(t_i))$.
We have %three
six  possibilities for $F_k=F_k(\xi_k^-(t_i), \xi_k^+(t_i))$:
% types of behavior of $w$.

\noindent
(a) $F_k$ is regular, does not have zero curvature, is non-interacting;\\
(b) $F_k$ is regular, does not have zero curvature, is interacting;\\
(c) $F_k$ is regular, has zero curvature, is non-interacting;\\
(d) $F_k$ is regular, has zero curvature, is interacting;\\
(e) $F_k$ is degenerate, does not have zero curvature, is non-interacting;\\
(f) $F_k$ is degenerate, has zero curvature, is non-interacting.

Cases (a) and (e) are solved in Lemma \ref{sanl2}, where corresponding
$\xi^\pm_k$ are constructed.

The construction of $\xi^\pm_k$ corresponding to (b), (d) is performed in Lemma
\ref{sanl3}. We stress that in all these cases  $\tau_k$, is
given by (\ref{sanr4}).

The definition of $\xi^\pm_k$ is simple if (c) or (f) holds, we just set
\begin{equation}\label{sanr17}
\xi^-_k(t) = \xi^-_k,\quad \xi^+_k(t) = \xi^+_k,\quad \tau_k(t)=0.
\end{equation}
We have to define $\Om$. By the very definition of zero-curvature
facets the intersection
$\partial J(\xi^+_k +\ep) \cap \partial J(\xi^-_k -\ep)$ is a
singleton $\{\al\}$ for any positive
$\ep< \min\{ \xi^-_{k+1} -\xi^+_k, \xi^-_k- \xi^+_{k-1}\}$. Moreover,
$\al\in\cA$, hence we set
\begin{equation}\label{oemga}
\Om(x,t) = \al,\quad\hbox{for } x\in [\xi^-_k(t), \xi^+_k(t)].
\end{equation}

Thus, we have specified evolution of $\xi^\pm_k$ for every configuration. In
all these cases the functions $\xi^\pm_k$, $k=1$, \ldots, $N_i$ are defined
on  maximal intervals $[t_i, t_i+T^\pm_k]$. The numbers $T^\pm_k$ are defined
as follows.

In (a) and (e) the positive number $T^+_k$ (resp.  $T^-_k$) is such that
$\xi_k^+(t) < \xi_{k+1}^-(t)$ (resp. $\xi_{k-1}^+(t) < \xi_{k}^-(t)$)
for $t< t_i+ T^+_k$ (resp.  $t< t_i+ T^-_k$), while equality occurs at
$t=T^+_k$ (resp. $t=T^-_k$), i.e. the facet begins to interact with
its neighbor. By Corollary \ref{sanc1} $T^\pm_k$ are finite.

If a group of interacting facets $F_l$, \ldots $F_{l+r}$ does not contain any
zero-curvature facet,
then by Proposition \ref{sanp1} it may not vanish and its maximal
existence time is defined as in (a) for $\xi_{l+r}$. Thus, at
$T^+_{r+l}$ the group begins to interact with another facet.
On the other hand, if this group of  interacting facets $F_l$, \ldots $F_{l+r}$
contains a zero-curvature facet, say $F_p$, then $T^+_p$ is defined as the
extinction time of $F_p$, i.e. $\xi_p^-(t)< \xi^+(t)$ for
$t\in[t_i,t_i+T^+_p)$, while $\xi_p^-(t_i+T^+_p)= \xi^+(t_i+T^+_p)$. Thus, the
number of facets drops by one.

Cases  (c) and (f) do not contribute to the definition of $t_{i+1}$, because
 (\ref{sanr17}) is valid for all $t\ge t_i$.

We have to define also $\Om(x,t)$. An attempt to do so reveals another
difficulty related to construction of $\xi_k^\pm$ starting from $t=0$.
Let us consider two interacting facets $F_k(\xi_k^-, \xi_k^+)$,
$F_{k'}(\xi_{k'}^-, \xi_{k'}^+)$,  where
\begin{equation}\label{sanr6}
[\xi_k^-, \xi_k^+]\subset \Xi_l(w_{0,s}),\quad
[\xi_{k'}^-, \xi_{k'}^+]\subset \Xi_r(w_{0,s}).
\end{equation}
%satisfying
%\begin{equation}\label{sanr7}
%[\xi_k^-, \xi_k^+]\cap [\xi_{k+1}^-, \xi_{k+1}^+] = \{\xi_k^+\}.
%\end{equation}
It is obvious that for any $s\in[\xi_k^-, \xi_k^+]$ and $s'\in [\xi_{k'}^-,
  \xi_{k'}^+]$  the intersection 
$$
\partial J(w_{0,s}(s)) \cap \partial J(w_{0,s}(s')) 
$$
is non-empty if and only if $|l-r|=1$. If the above intersection is non-empty,
we can construct the desired 
$\Om(x,t)$. On the other hand, if this intersection is void, then we have no
chance to construct a $W^{1}_{1}$ section of $\partial J(w_s)$.

Let us suppose then that (\ref{sanr6})  holds and $|l-r|=p+1$,
$p>0$. Let us suppose for simplicity that $l<r$. Thus, a single point
$\xi$ is a connected component of $\Xi_j(w_{0,s})$, $j=l, l+1,\ldots,r$, i.e.
$$
\xi = \xi_j^-= \xi_j^+, \quad j=l, l+1,\ldots,r.
$$
In other words, we have a number of degenerate, interacting facets at
$\xi$. The system of ODE's (\ref{sanr12}) is singular. The  problem of
evolution of interacting degenerate facets shall be dealt with below
in Part B of the proof. It occurs only at $t=0$.

Finally, we check that $w(x,t)$ and $\Om(x,t)$ fulfill the conditions
postulated in the definition of the weak solution. They satisfy the
equation
\begin{equation}\label{bullet}
w_t (x,t)= \Om_x (x,t)
\end{equation}
and the initial and boundary conditions are satisfied.
Integral identity in Definition \ref{defi} follows.

\bigskip
{\sc Part B. }
After finishing part A, i.e. the case of data  satisfying
(\ref{sanr13}), we
consider the interaction of facets during creation,
i.e. (\ref{sanr13}) is no longer valid. We have two cases to consider:\\
(g) $F_k$ is  degenerate, with nonzero curvature and interacting;\\
(h) $F_k$ is  degenerate, with zero curvature and interacting.\\

We begin with (g). Let us suppose that $w_0$ violates
(\ref{sanr13})  at some $\xi$.
Thus, we are dealing with the situation when one sided derivatives of
$w_0$ differ at $\xi$, i.e.,
$$
w_{0,s}^-(\xi) < \al_k < w_{0,s}^+(\xi)
$$
for some $a_k\in\cA$. It may as well happen that the reverse inequalities occur,
however for the sake of definiteness we shall stick to the above choice.

We shall construct two functions $w_\ep$, $w^\ep$ such that their derivatives
belong to $\JR$, $w_\ep(x) < w_0(x) < w^\ep(x)$ and
%$\ep>0$ such that
\begin{equation}\label{rnuz2}
|w_\ep(x) - w_0(x)|,\quad | w^\ep(x) - w_0(x)| <\ep. 
%,\qquad |x-\xi^+_k|\ge \ep,
\end{equation}
%$$ w_\ep(x) < w(x) < w^\ep(x),\qquad |x-\xi^+_k|< \ep.$$
%Moreover we require that $[\xi^-_{k-1,\ep},\xi]$,
%$[\xi,\xi^+_{k,\ep}]$ (resp. $[\xi^{-,\ep}_{k-1},\xi]$,
%$[\xi,\xi^{+,\ep}_{k}]$) are
%connected components of $\Xi(w_{\ep,s})$ (resp.  $\Xi(w^\ep_s)$) and
%$$
%|\xi^-_{k-1,\ep}-\xi|,\quad
%|\xi-\xi^+_{k,\ep}|,\quad |\xi^{-,\ep}_{k-1}-\xi|,\quad
%|\xi-\xi^{+,\ep}_{k}|\le \ep/4.
%$$
We set
$$
w^\ep(x) = \max \{w_0(x), l_k(x,\xi, w_0(\xi)+\de)\},
$$
where $\de>0$ is so chosen to guarantee (\ref{rnuz2}).
We also define
$$
w_\ep(x) = \max \{w_0(x)-\ep, l_k(x,\xi, w_0(\xi)+\de)\},
$$
where $\de\in(0,\ep)$ is arbitrary. Of course,  (\ref{rnuz2}) holds.

If the newly constructed $w^\ep$ and $w_\ep$ do not satisfy (\ref{sanr13}), we
repeat the above process until they do. Subsequently,
we apply the results of Part A to $w^\ep$ and $w_\ep$. We deduce from that
existence of interacting facets at $\xi$. By the comparison principle,
the non-zero interacting facets exist for $t>t_0$.

Finally, we study (h). We notice that at such an instance $F_k$ cannot interact
with two neighboring facets, because this would mean that  $F_{k-1}$ and
$F_{k+1}$ lay on the same line, that is, $F_k$ is their common end point. Thus,
the three facets $F_{k-1}$, $F_k$ and $F_{k+1}$ form a single facet $\tilde
F_k$ with the pre-image $[\xi^-_{k-1}, \xi^+_{k+1}]$. On the other hand it may
happen that $F_k$ is a degenerate, zero curvature facet interacting
with just one neighbor, say $F_{k+1}$. 
Since $F_k$ is degenerate, i.e., $\xi_k^+=\xi_k^- =:\xi_k$, due to its
interaction with $F_{k+1}$ we have $\xi^-_{k+1} = \xi_k$. Moreover, 
%This happens when 
$w_{0,s}^-(\xi_k) =\al_{k-1}$ %and <\al_k$
and $w_{0,s}^+(\xi_{k}) =\al_{k}$ where $\al_{k-1}$, $\al_k\in\cA$ and 
we may assume that $\al_{k-1} < \al_k$, (the other case is handled similarly)
and $w_s(\xi^+_{k-1},\xi_k)\subset
(\al_{k-2},\al_{k-1})$. A similar situation occurs when $F_k$ interacts with
$F_{k-1}$.

In order to determine the evolution of the system we have to take into account
if $F_{k+1}$ has zero-curvature or not. In the former case $\dot\tau_{k+1}
=0$, hence we set $\tau_k\equiv0$. In latter case we have $\dot\tau_{k+1}>0$
(it may not occur $\dot\tau_{k+1}<0$). Thus, $F_k$ disappears instantly. As a
result, we agree to disregard $F_k$ and diminish $N_0$ by 1.

\bigskip
{\sc Part C.}
We have to deal with the points outside of \, $\Xi(w_s(\cdot,t))\equiv
\bigcup_{i=1}^{N_k}[\xi^-_{k},\xi^+_k]$. By the definition of
$\Xi(w_s(\cdot,t))$, its complement is open
$$
[0,2\pi]\setminus \Xi(w_s(\cdot,t)) =\bigcup_{l=1}^{N_i} (\xi^+_k,\xi^-_{k+1}).
$$
where $(\xi^+_{N_i},\xi^-_{{N_i}+1})$ should be understood as
$(\xi^+_{N_i},2\pi]\cup [0, \xi^-_1)$, (with the understanding that
  $0\le\xi^\pm_k \le 2\pi$, $k=1,\ldots, N_i$).
Using again the definition of $\Xi$, we come to the conclusion that, if $x$
belongs to any of the intervals $ (\xi^+_k,\xi^-_{k+1})$, then either
$w_s(x,t)$ exists or $w_s^+ (x,t)\neq w_s^-(x,t)$. In either case, the set
$\partial w(x,t)$ (see (\ref{depaw})) does not intersect $\cA$. Since
$\partial w(x,t) $ is an
interval, we deduce that there exists $\al_k\in\cA$ such that
\begin{equation}\label{rnalf}
\partial w(x,t) \subset (\al_k, \al_{k+1}).
\end{equation}
We have to make sure that the choice of $\al_k$, in the formula above, depends
only on the interval $(\xi^+_k,\xi^-_{k+1})$, but it is independent
from a specific point $x\in (\xi^+_k,\xi^-_{k+1})$. Indeed, by the
definition of the $\JR$ class $\partial w =
M-f$ or $\partial w = f -M$ , where $f$ is a continuous increasing function
and $M$ a maximal monotone operator. Thus, the images
$f(\xi^+_k,\xi^-_{k+1})$ and $M(\xi^+_k,\xi^-_{k+1})$ are connected
intervals, so is the image $\partial w(\xi^+_k,\xi^-_{k+1})$, which is
disjoint from $\cA$. Our claim follows.

As a result, our definition of $w(x,t)$ for $x \not \in \Xi(w_s(\cdot,t))$ is
as follows,
$$
w(x,t) = w(x,t_k) \quad\hbox{and}\quad\Om(x,t) = \frac{dJ}{d\vfi}(w_s(y,t_k))
\qquad\hbox{for }x \in (\xi_i^+, \xi_{i+1}^-).
$$
where $y \in (\xi_i^+, \xi_{i+1}^-)$ is any differentiability
point of $w(\cdot,t_k)$.

\bigskip
{\sc Part D.}
We have to define $t_{k+1}$. We do this inductively. Once $t_k$ is given, we
set
$$
t_{k+1} = t_k + \min\{\min_i T^+_i,\min_i T^-_i\}.
$$
Thus at $t_{k+1}$ two facets begin to interact, due to the shrinkage of
$[\xi_i^+, \xi_{i+1}^-]$ to a point or due to the disappearance of a
facet. By Proposition \ref{sanp1}, we know that only zero-curvature
facets may disappear. We set 
$$
N_{i+1} = N_i -m,
$$
where $m$ is the number of removed degenerate, interacting, zero-curvature
facets at $t=t_{i+1}$.

The last thing to show is the estimate $\| w_s(\cdot,t)\|_{\JR} \le\|
w_s(\cdot,s))\|_{\JR}$, whenever $t>s$. By the construction above, the
number of connected components of $\Xi(w_s(\cdot,t))$ drops at time
instances  $t_k$, $k=1,$\ldots, $M_N$, hence $K(w_s(\cdot,t))\le
K(w_s(\cdot,s))$, whenever $s\le t$. It remains to show that
$\|w_s(\cdot,t)\|_{TV(S)}\le 
\|w_s(\cdot,s)\|_{TV(S)}$, where we denoted by $ \|f\|_{TV(E)}$ the total
variation of function $f$ over set $E$.

We first consider the case $t>s$ such that $\Xi(w_s(\cdot,t))\ne S$,
we know that we always have $\Xi(w_s(\cdot,t))\supset
\Xi(w_s(\cdot,s))$ for $s<t$. By the general properties of the total
variation, we notice that
$$
\|w_s(\cdot,t)\|_{TV(S)} =\|w_s(\cdot,t)\|_{TV(\Xi(t))} +
\|w_s(\cdot,t)\|_{TV(S\setminus\Xi(t))} ,
$$
where we wrote $\Xi(\si)$ for $\Xi(w_s(\cdot,\si)$. Now, by the
definition of $w(x,t)$, we notice that
$$
\|w_s(\cdot,t)\|_{TV(S\setminus\Xi(t))} =
\|w_s(\cdot,s)\|_{TV(S\setminus\Xi(t))} \le
\|w_s(\cdot,s)\|_{TV(S\setminus\Xi(s))}.
$$
We turn our attention to $\|w_s(\cdot,t)\|_{TV(\Xi(t))}$. On the
intervals forming  $\Xi(t)$ function $w_s(\cdot,t)$ is piecewise
constant. The jumps occur at the endpoint of these intervals. They are
no bigger and no more numerous than the jumps of
$w_s(\cdot,s)$. Thus our claim follows in the considered case of
$t$. In fact, the case of $t$ such that $\Xi(t)=S$ is not much different.
Finally, we can see that $w_s$ is a difference of two monotone
functions and one of them is continuous, the other one a maximal
monotone operator.

Our theorem is proved.
\qed

We close this subsection with a formula, which might be called
``morphing a circle into a square''.

\noindent{\bf Example.} Let us suppose that $\phi_0(s) = s$ or
$w_0(s)=\frac12 s^2$. Due to the high symmetry of the problem, it is
sufficient to consider just formation of one facet. Then, $w(x,t)$, the
unique solution to (\ref{nasz}), is given by the formula,
\begin{equation*}
w(x,t)  =\left\{
\begin{array}{ll}
\frac12 s^2 &  s\in [0,\xi^-_1(t)] \cup[\xi^+_1(t),\frac\pi 2],\\
\frac\pi4 s -\frac{\pi^2}{32} +\tau_1(t) & s\in [\xi^-_1(t), \xi^+_1(t)].
\end{array}\right.
\end{equation*}
Here, $\xi^\pm_1 = \frac\pi4  \pm \sqrt{2\tau_1}$ and $
\tau_1 = \left(\frac{\sqrt 2\pi }{12}t\right)^{2/3}$.
Let us note that at $T_1=\pi^2/2^6$ we have $\Om^+-\Om^- =
\xi^+-\xi^-$, so for later times $\dot\tau =1$.

\bigskip
We can make this observation more general.

\begin{proposition}\quad{\sl Let us suppose that the assumptions of Theorem
    \ref{sant1} are satisfied. Then, there exist $T_{fa}$, such that
  if $t>T_{fa}$, then $w(\cdot,t)$ is fully faceted, i.e. $w_s(\cdot,t)$ is
  piece-wise linear. More precisely, for $\Xi(w_s(\cdot,t))\subset[0,2\pi)$ for
  $t<T_{fa}$ and $\Xi(w_s(\cdot,t))  =[0,2\pi)$ for
  $t\ge T_{fa}$.
}
\end{proposition}

\bigskip\noindent{\it Proof.} Let us consider $w_0$. It is fully
faceted or not. If 
not, then by the proof of Theorem \ref{sant1}, we  deduce that after at some
$t_{i_0}$ we have $\Xi(t_{i_0})=[0,2\pi)$ and our claim follows.\qed

\subsection{Convexification}
We show that after some depending upon the initial data, the solution becomes
such that $w_s = \vfi$ is monotone decreasing or increasing. We shall call this
process by convexification.

\begin{proposition} \quad
{\sl  Let us suppose that the assumptions of Theorem
    \ref{sant1} are satisfied. Then, there exist $T_{cx}$, such that
  if $t\ge T_{cx}$, then $w_s(\cdot,t)$ is monotone, while this is not true for
    $t<T_{cx}$}.
\end{proposition}

\noindent{\it Proof.} If $w_{0,s}$ is monotone, then we are done.
Otherwise, let us
suppose that $t_j$ is the largest time such that at $t_j$ a zero curvature
facet disappears. Since the zero-curvature facets cannot persist because their
endpoints necessarily move, it follows that $T_{cx}= t_j$ has the desired
properties. \qed

\bigskip\noindent{\bf Remark.} All possibilities can be realized
$T_{cx}> T_{fa}$ as well as $T_{cx}< T_{fa}$.

\subsection{Asymptotic behavior of facets}% interactions}

Here, we consider the last stage of evolution, when $t\ge t_M$ and
$N_M=4$. In this case, it is sufficient to specify only $\xi_k^+$,
$k=1,2,3,4$. Furthermore, the system for interacting facets,
(\ref{sanr12}) takes the form,
\begin{eqnarray}\label{sanr16}
&& \al_1\dot \xi_1 -  \al_2\dot \xi_2 = \dot\tau_1 -\dot\tau_2\nonumber\\
&& \al_2\dot \xi_2 -  \al_3\dot \xi_3 = \dot\tau_2 -\dot\tau_3\nonumber\\
&& \al_3\dot \xi_3 -  \al_4\dot \xi_4 = \dot\tau_3 -\dot\tau_4\\
&& \al_4\dot \xi_4 -  \al_1\dot \xi_1 = \dot\tau_4 -\dot\tau_1\nonumber\\
&& \xi_k(t_M) = \xi_k,\quad k=1,2,3,4.\nonumber
\end{eqnarray}
We notice that the stationary points of (\ref{sanr16}) are such that
$\dot\tau_1=\ldots=\dot\tau_4$. This occurs if and only if
$\Om^+_k-\Om^-_k = \xi_{k} -\xi_{k-1},\quad k=1,2,3,4,$
where by $\xi_0$ we understand $\xi_4$. Moreover, due to our
assumptions on $J$ we  have $\Om^+_k-\Om^-_k =\De\Om$, $k=1,2,3,4$.

Additionally, system  (\ref{sanr16}) possesses a Liapunov
functional. Namely, let us write
$$
F(\vec\xi) = \sum_{k=1}^4 \ln(\xi_k-\xi_{k-1})^{\De\Om},
$$
with the understanding of $\xi_0$ as above. By direct calculation, we
check that
%We also define $D$ to be a diagonal matrix with the entries
%$(\al_k-\al_{k-1})^{-1}$. Then, we notice
$$
\frac d{dt} F(\vec\xi) = \nabla_\xi F\cdot \frac d{dt} \vec\xi
%= - \nabla_\xi F\cdot A \nabla_\xi F
< 0.
$$
This derivative vanishes if and only if  $\xi$ is the only equilibrium
point. Thus, we have a complete picture of the asymptotic behavior of
$\La$.

\begin{theorem}{\sl Let us assume that
$\vfi_0\in\JR$ and $w$ is the corresponding  unique
    solution to (\ref{nasz}). Then, there exists $T_1$,
$\max\{T_{cx},T_{fa}\}\le T_1\le\infty$ with the following property:\\
(a) If $T_1<\infty$, then
$\xi_l(t) = -\frac{3\pi}4+ \frac\pi2 l +\al,$ for some $\al\ge0$,
    $l=0,$ \ldots, 3, and
$t\ge T_1$, in other words, $w$ is the minimal solution for $t>T_1$;\\
(b) If $T_1=\infty$, then $\lim_{t\to\infty}\xi_l(t) = -\frac{3\pi}4+
    \frac\pi2 l
    +\al,$ $l=0,$ \ldots, 3 for some $\al\ge0$.}
\end{theorem}\qed

\subsection{Proof of Theorem \ref{thmreg}}

In the course of proof of Theorem \ref{sant1}, we exhibited a quite
explicit construction of the weak solution with such initial data that
$\vfi_0\in \JR$. Now, we have to show that is has all the postulated properties
of the almost classical solution. We have already noticed that $w_s =
\La_s + s$ belongs to the \JR  \ class, furthermore
$\|w_s(\cdot,t)\|_{\JR}\le \|w_s(\cdot,0)\|_{\JR}$. The key point,
however, is to realize that
\begin{equation}\label{jeje}
\Om =\partial J \bar\circ\partial w,
\end{equation}
where $\partial w$ is the multivalued map whose section is $w_s$. We
defined  $\partial w$ in (\ref{depaw}).
Checking that (\ref{jeje}) indeed holds requires recalling the steps of
construction of $\Om$, we will do this below. Finally,
after we set $N=\{0, t_1,\ldots, t_M\}$ we see that
$$
\La_t = \frac{\partial}{\partial s} {\partial}J\bar\circ(\La_s+s),
$$
holds for all $t\in (0,+\infty)\setminus N$ in the $L^1$ sense, more precisely
it holds pointwise except $x\in
[0,2\pi)\setminus\{\xi^\pm_i:\ i=1, \ldots, N_k\}$. Indeed, the
definitions (\ref{sanr4}), (\ref{uta}),
(\ref{sanr17}) of $\tau_k(t)$ were such that
$\frac d{dt}\tau_k(t) = \frac{\partial\Om}{\partial s}$. Moreover,
$\frac{\partial\La}{\partial t} = \frac d{dt}\tau_k(t)$, see Lemma
\ref{sanl2}, Lemma \ref{sanl3} and eq. (\ref{bullet}). We recall
that by definition functions $\tau_k(\cdot)$ are continuous on
$[t_i,t_{i+1}]$ and differentiable in $(t_i,t_{i+1})$. Moreover, the
right derivative of $\tau_k(t)$ is  well-defined for all $t$, except
possibly  $t=t_0$. Hence, $\frac{\partial\La}{\partial t}$ is 
defined everywhere, except the points $t_i$, $i=0,\ldots, M$, but the
right time derivative $\frac{\partial\La^+}{\partial t}$ is defined
for all $t>0$.

We will check below that $\Om$, constructed in the course of proof of Theorem
\ref{sant1}, coincides with $\partial J\bar\circ \partial w$,  --- see
(\ref{sanr5}), (\ref{omega}), (\ref{oemga}), where $w_s(s,t)
= \La_s(s,t)+s$.
In order to see that
we examine the
steps of the construction of $\Om$ and compare it with the definition of the
composition $\bar\circ$. Let us fix $t\in(t_k,t_{k+1})$, at the end we will
consider $t=t_{k+1}$, then we compose $\partial w(\cdot,t):[0,2\pi] \to [a,b]$
with $\partial J:\bR\to \bR$. We have to identify the sets $\cD_s$, $\cD_f$
and $\cD_r$ appearing in the Definition \ref{def3}. For our  choice of $t$ we
have
$$
\cD_s (t) = \{ s\in [0,2\pi]: w^+_x (s,t) \neq w^-_x(s,t)\}.
$$
In particular, $\cD_s (t)$ contains all points $\xi^\pm_i(t)$,
$i=1,\ldots,N_k$. We can see that
$$
\cD_f (t) = \bigcup_{i=1}^{N_k} (\xi_i^-(t),\xi_i^+(t)),
$$
i.e., it is the sum of interiors of intervals contained in $\Xi(w_s(\cdot,t))$.
Finally, by the definition
$$
\cD_r  (t)= [0,2\pi] \setminus (\cD_s (t) \cup \cD_f (t)).
$$
We shall consider these cases separately.

$1^o$ case $\cD_r $. If $s\in  \cD_r  (t)$, then $w$ is differentiable
at $s$ and $w_s(s,t)\not\in \cA$. Thus, by (\ref{de2}) $\partial
J\bar\circ\partial w(s,t)= \frac {dJ}{d\vfi}(w_s(s,t))$. We notice that
$\cD_r (t) \subset [0,2\pi] \setminus \Xi(w_s(\cdot, t))$, hence by
Part C of the proof of Theorem \ref{sant1} we immediately see that  $\partial
J\bar\circ\partial w(s,t)$ equals $\Om(s,t)$ on $\cD_r $.

$2^o$ case $\cD_f$. By its definition  $\cD_f(t)$ is the sum of
interiors of pre-images of facets, as noticed above. Moreover, on each
interval $(\xi_i^-(t),\xi_i^+(t))$, the set $\partial w(x,t)$ is a
singleton equal to $\{\alpha_k\}\subset \cA$. Then, the cases of the Definition
\ref{def3}, see formulas (\ref{de3})--(\ref{de6}) have their
counterparts in the formulas (\ref{sanr5}), (\ref{omega}) and (\ref{oemga}).

$3^o$ case $\cD_s$. We notice that, if $t>0$, then the set
$\Xi(w_s(\cdot,t))$ has no component, which is a singleton. Thus,
if $s \in  \cD_s(t)$, then the set $\partial w(s,t)$ does not
intersect $\cA$. As a result, formula (\ref{de7}) for the composition
yields a singleton, because on the RHS of (\ref{de7})  the limit of
constant functions are taken. This in agreement with the discussion of
Part C.

Finally we have to deal with the case $t=t_{k+1}$. On one hand
$\Om(\cdot,t_{k+1})$ is defined by the left time continuity of $\Om$, on the
other hand we have to check that $\Om =\partial  J\bar\circ\partial
w$.

By the very definition of $t_{k+1}$ (see Part D of the proof of Theorem
\ref{sant1}), at this time instant a zero-curvature  curvature facet
disappears or two facets  begin to interact or merge, i.e.,
$$
\lim_{t\to t_{k+1}^-} \xi_i^+ (t) = a = \lim_{t\to t_{k+1}^-}
\xi_{i+1}^- (t).
$$
We have then two possibilities, either
$a = \xi_j^+ (t_{k+1}) = \xi_{j+1}^- (t_{k+1})$ or
$a\in (\xi_{j}^- (t_{k+1}), \xi_j^+ (t_{k+1}))$ where this interval is
a connected component of $\Xi(w_s(\cdot, t_{k+1}))$. Once we realize
this, it is clear that
$\Om(\cdot, t_{k+1})= \partial J\bar\circ\partial w(\cdot, t_{k+1})$.
\qed
%##############################3

%%% --- s 7 ---
\section{Appendix}

\subsection{Motivation of equation (\ref{rn2})}

Here,
we consider closed curves, we view them as graphs over a smooth, convex
reference closed curve $\cM$.  We do not make here any attempt to consider
non-smooth reference curves, which is reasonable because this would
add up difficulties while not giving advantages.

Let us suppose that $x_0(s)$ is an arc-length parameterization of $\cM$ and
$\te_t(s)$, $\te_n(s)$ are unit tangent and normal vectors,
respectively, such that $(\te_t(s),\te_t(s))$ is positively oriented. Then all
points in a neighborhood of $\cM$ can be uniquely written as $x=
x_0(s) + \te_n \La$, as a result we can parameterize our curve $\Ga(t)$
as
$$
x(s,t) = x_0(s) + \te_n(s) \La(s,t).
$$
Since $\cM$ is convex we may write $\te_n$ uniquely as
$\te_n(\vfi(s)) = (\cos \vfi(s), \sin \vfi(s))$, where $\vfi$ is the
measure of the angle between the $x_1$ axis and $\te_n$. Moreover,
$$
\frac d{ds} \te_n(\vfi(s)) = - \te_t(\vfi(s)) \frac {d\vfi}{ds}=
 - \kappa\te_t(\vfi(s)).
$$

We note
$$
\frac{\partial x}{\partial s}(s,t) =
\te_t (1 - \ka\La) +
\te_n \La_s,
$$
because $|\dot x_0(s)| = 1$.
With this formula at hand, we can write the expression for the tangent
and normal to $\Ga(t)$, they are
${\bf \tau} = \frac 1W(
\te_t(1 - \ka \La) +\te_n \La_s)$,
$ \bn =  \frac 1W( - \La_s\te_t +(1 - \ka \La)\te_n )$,
where
$W^2 = (1 - \ka \La)^2 + \La_s^2$. Hence, the LHS of
(\ref{rn0}) takes the form
$$
\be V = \be \frac{dx}{dt}\cdot\bn =
\frac 1W (1 - \ka \La) \La_t.
$$
The RHS of (\ref{rn0}) is $\hbox{div}_S\nabla_\xi\ga(\xi)|_{\xi=\bn}$.
In our paper \cite{MuchaRybka}, we have shown that  it is equal to
$$
\ka = \frac
d{ds}\left(\frac{\partial}{\partial\vfi}I_\theta(\vfi)\right).
$$
We defined $I_\theta(\vfi)$ as follows,
$I_\vartheta(\vfi) = \bar\gamma(\bn(\vfi))+
\int_{\vartheta}^{\vfi}d\psi\int_{\vartheta}^\psi \bar\gamma(\bn(t))dt.$
We noted that this function is convex iff the stored energy function
$\bar\gamma$ is convex. However, in general $I_\vartheta$ does not enjoy
higher regularity properties. It is not differentiable at angles
corresponding to the normals to the Wulff shape.

Finally, equation
(\ref{rn0}) takes the form
\begin{equation}\label{rn1}
\be \bn\cdot\te_n \La_t =
\frac d{ds}\left(\frac{\partial}{\partial\al}I_\theta(\al)\right),
\end{equation}
where $\al$ is the measure of the angle between
the $x_1$ axis and $\bn$.

One may study evolution of convex curves defined by their angle
parameterization. We notice
$\al = \vfi + \psi,$
where $\psi$ is the measure of the angle between ${\bf \tau}$ and $\te_t$.
We notice that
${\bf \tau }\cdot \te_n = \sin\psi  = \frac{\La_s}W$,
${\bf \tau } \cdot \te_n = \cos\psi  = \frac{1-\La}W.$
Thus, we see that $\psi = Arg({\bf \tau} \cdot \te_n + i {\bf \tau} \cdot \te_n)$,
in fact we have
$\psi = \arctan \left(\frac{\La_s}{1-\La}\right).$
Thus, (\ref{rn0}) takes the form
$$
\be \bn\cdot\te_n \La_t =
\frac d{ds} \left(
\frac{\partial}{\partial\phi}I_\theta
\left(\vfi+\arctan\left(\frac{\La_s}{1-\La}\right)\right)
\right).
$$
This equation is rather involved, we prefer to simplify it by dropping
the terms which at this stage we deem not important, thus we come to
(\ref{rn2}).

\subsection{Other choices of function $J$}

We may also consider any properly chosen piecewise linear, convex
$J$,
\begin{equation}\label{jotGL}
J_l(\vfi) = \sum_{i=1}^N b_i |\vfi - \al_i|.
\end{equation}
We require that $N\ge 4$, $b_i>0$ and
$\al_0<\al_1<\ldots<\al_N<\al_0+2\pi$, we will write
$S=[\al_0,\al_0+2\pi)$. In order to stick to geometrically relevant
data, we also impose the condition that $\sum_{i=1}^N b_i =\pi$, which
guarantees that $\partial I(S)$ is an interval of length $2\pi$. In
addition, we assume that
the following function yields an
angle parameterization of closed curve, which encompasses a convex
region. Namely, we set
\begin{equation}\label{Omka}
\Om_j = \sum_{i=1}^j b_i - \sum_{i=j+1}^N b_i,\qquad j=0,\ldots, N,
\end{equation}
with the convention that the summation over an empty set of parameters
yields zero. Then we define
$\Phi:[\al_0,\al_0+2\pi)\to \bR$ by the formula
\begin{equation}\label{wilk}
\Phi(s) = \sum_{i=0}^N
\Om_i\chi_{[\al_{i},\al_{i+1})},
\end{equation}
(with the convention $\al_{N+1}=\al_0+2\pi$) is an
angle parameterization of closed curve. We notice that our
assumptions imply that $\Om_0 +2\pi = \Om_N$.

The analysis of behavior of solutions presented in Section \ref{sas} is 
valid also for $J$ given by (\ref{jot}) and $J_l$, however the actual
calculations for $J_l$ are more lengthy. In addition we
may show existence of weak solution for a general, piecewise smooth, convex
$J$, but in this case we cannot offer detailed analysis of solutions,
yet.

\bigskip\noindent{\bf Acknowledgment.} The present work has been
partially supported by Polish KBN   grant No. 1~P03A~037~28. Both
authors thank Universit\'e de Paris-Sud XI, where a  part of
the research for this paper was performed, for its hospitality.


\begin{thebibliography}{LSU2}
\bibitem[L]{luckhaus} S.Luckhaus,  Solutions for the two-phase Stefan
problem with the Gibbs--Thomson law for the melting temperature,
{\it Euro. J. Appl. Math.,} {\bf 1}, (1990),  101--111.
\end{thebibliography}
\end{document}